\documentclass[12pt]{article}

\RequirePackage{a4}
\usepackage{epsfig}
\usepackage[utf8x]{inputenc}
\usepackage{verbatim}
\usepackage{lineno}

\usepackage{hyperref}
\usepackage{url}
\usepackage{breakurl}

\begin{document}

\title{ Talking about Probability, Inference and Decisions\footnote{\,To
    be published in {\em Progetto Alice,
      Rivista di matematica e didattica}.} \\
  --\  Part 1: The Witches of Bayes\  -- \\
  (In Italian)
}

\author{Giulio D'Agostini\footnote{\,giulio.dagostini@roma1.infn.it,
 \url{http://www.roma1.infn.it/\~dagos} (contact author)} \\
Dipartimento di Fisica, Universit\`a di Roma ``La Sapienza'' \\ \mbox{} \\ 
Noemi Cifani\footnote{\,noemi.cifani@gmail.com}\,
and Alba Gilardi\footnote{\,alba.gilardi@alice.it} \\
Istituto Comprensivo Giuliano Giorgi di Montorio Romano (RM)\\
Scuola Media Plesso di Monteflavio (RM)
}

\date{}

\maketitle

\begin{abstract}
In October 2017 the Italian National Institute
of Statistics (ISTAT), Italy’s body for official statistics,
has published the book of fairy tales {\em Le streghe di Bayes}
({\em The witches of Bayes})~\cite{Streghe-ISTAT},
written by ISTAT staff members with the commendable aim of
introducing statistical and probabilistic reasoning
to children. In this paper the fairy tale which gives the name
to the book is analyzed in a dialog between three teachers with different
background and expertise. The outcomes are definitively discouraging,
especially when the story is compared to the appendix of the
book, in which the teaching power of every story is indeed
explained\footnote{\,For an English translation of the appendix,
  together with a synopsis of the story, 
  see \cite{ISBA-Bulletin}.}
(as a matter of fact, without the appendix the fairy tale of the witches
seemed to be written with the purpose of make the `Bayesians', 
meant as the villagers from `Bayes', ridiculous).
In fact the fairy tale of the witches
 does not contain any Bayesian reasoning, the suggested decision
strategy is simply wrong and the story does not even
seem to be easily modifiable (besides the trivial correction of the
decision strategy) in order to make it usable as a teaching tool.
As it happens in real dialogues, besides the
fairy tale in question, the dialogue touches
several issues somehow related to the story and concerning
probability, inference, prediction and decision making.
The present paper is an indirect response to the invitation by the
ISBA bulletin to comment on the fairy tale~\cite{ISBA-Bulletin}.
\end{abstract}

\newpage
\begin{center}
  {\Large Conversando su probabilità, \\ inferenza e decisioni } \\
  \mbox{} \\
       {\huge --\,  Parte 1: Le streghe di Bayes\, --}
       \mbox{} \\  \mbox{} \\  \mbox{}  \\ \mbox{}  \\\mbox{} \\
            {\bf Sommario} 
\end{center}
\mbox{}\hspace{5mm}A ottobre del 2017 l'ISTAT ha pubblicato
{\em Le streghe di Bayes}~\cite{Streghe-ISTAT}, un
libro di fiabe scritto da personale dell'Istituto
con l'encomiabile intento di introdurre ai bambini
il ragionamento statistico e probabilistico.
La fiaba che dà il nome al libro è qui analiz\-zata in forma
di dialogo fra tre insegnanti con diverso background e esperienza.
Il risultato è decisamente sconfortante, soprattutto quando
si confronta la fiaba con l'appendice, in cui vengono illustrate
le potenzialità didattiche di ciascuna storia 
(senza l'appendice la fiaba delle streghe sembrerebbe essere
stata scritta con l'intento di mettere in ridicolo i `bayesiani', 
intesi come gli abitanti del villaggio di fantasia `Bayes'). Infatti
la fiaba delle streghe non contiene alcun ragionamento bayesiano, la strategia
decisionale proposta è semplicemente errata e la storia non sembra
nemmeno facilmente modificabile (a parte la banale correzione
della strategia decisionale) in modo tale da renderla usabile
per l'insegnamento.
A parte la fiaba in questione, come avviene nei dialoghi reali,
il dialogo tocca 
diversi aspetti riguardanti probabilità, inferenza, previsione
e decisione collegati in qualche modo alla storia.
Questo scritto è una risposta
indiretta all'invito del Bollettino della Società Internazionale
di Analisi Bayesiana di riportare
commenti sulla fiaba delle streghe\,\cite{ISBA-Bulletin}. 
\mbox{}\\ 
\mbox{}\\ 
\section*{}
\begin{description}
\setlength\itemsep{0.5mm}
\item[G.] Allora, dicevate, com'è che siete arrivate a questo libricino?
\item[A.] Noemi, appena arrivata nella nostra scuola,
  si è mostrata entusiasta del progetto di lettura di Istituto e
  mi ha chiesto se potevamo trovare un argomento che coinvolgesse
  anche la Matematica, che lei insegna.
  Così, sbirciando nella rete, abbiamo visto che si parlava di questo
  libro~\cite{Streghe-ISTAT}. L'idea del mix di fiabe e statistica ci è piaciuta subito. Tutti parlano di statistica di qua e statistica di là, della probabilità di questo e della probabilità di quello\ldots 
\item[N.] \ldots e nell'ignoranza ci sono quelli che vanno appresso ai numeri ritardatari.
\item[A.] Lo fanno perché, da quanto ho capito, c'è una legge della statistica che dice che tutti i numeri devono uscire con la stessa frequenza e quindi il numero in ritardo si deve mettere in regola, per così dire. 
\item[G.] Ah, è così che ti sembra?
\item[A.] Non è che sembra a me. Lo dicono tutti. E dicono che questa legge discende dalla teoria matematica della probabilità, insomma una cosa seria, a sentire loro. E i numeri in ritardo sono pure evidenziati sui giornali e nei posti dove si gioca. Ma a me, lasciamelo dire, mi sembra una grande sciocchezza, perché i numeri stanno dentro quella specie di gabbia che gira, e non capisco proprio come possa essere che quello in ritardo sgomiti per farsi strada verso la feritoia.
\item[G.] Divertente questa dei numeri ritardatari che sgomitano. Lo sai che ci sono anche altri numeri che fremono per uscire? Ma non per rispettare la legge, bensì per violarla, perché è nella loro natura, come direbbe Aristotele: i cosiddetti numeri caldi, o così almeno venivano chiamati quelli che erano usciti molto più spesso degli altri. E anche questi venivano segnalati, e la gente ci puntava su, preferendoli, insieme ai ritardatari, a quelli che invece erano più o meno in regola sulla `tabella di marcia'.
\item[A.] Ma perché parli al passato? Ora non sono più popolari?
\item[G.] Non so. Mi ricordo di cose viste decine di anni fa. In particolare c'era uno sull'autobus che mi faceva vedere non so quale settimanale specializzato che riportava tutte queste informazioni. Sto parlando di quando le estrazioni avvenivano settimanalmente. Ma credo che ci saranno ancora pubblicazioni del genere, oggigiorno più probabilmente su appositi siti web.
\item[A.] Però capite bene che quando i numeri ritardatari  --  non so se pure quelli che tu chiami `caldi', la prossima volta ci farò caso  --  sono bene evidenziati in bar e tabaccherie, sono modi per indurre la gente a giocare. E poi si parla di ludopatia! Ma lasciamo stare\ldots 
\item[G.] Non per niente qualcuno definiva il lotto
  una tassa sulla stupidità umana\footnote{\,La definizione
    è da taluni attribuita
    al matematico Bruno de Finetti, come si può
    verificare da una rapida ricerca su internet.
    Ma non sembra il solo. Una espressione del genere
    è attribuita anche a Garibaldi, ma l'argomento è
    al di là dei nostri interessi e dello scopo
    di questo scritto.}
  e molti stati dell'Italia prima dell'unificazione lo avevano,
  almeno temporaneamente, abolito.
  Ma poi, considerando che era una buona fonte di gettito per l'erario, si sono scrollati di dosso le remore basate su principi morali.
\item[A.] Sembra che l'abolizione del gioco del lotto
  sia stata una delle prime misure
  prese da Garibaldi dopo il suo ingresso a Napoli.
\item[G.] Questa non la sapevo. 
\item[N.] E il popolo? Ho difficoltà a immaginare i napoletani
  privati del lotto.
\item[A.] Il popolo\ldots\ già, il popolo\ldots\ Come si dice, il popolo vuole
  essere ingannato, ma ha anche il diritto di sognare\ldots\ Basta 
  che uno sia cosciente di quello che sta facendo e che il
  gioco non diventi una mania. Insomma, come facciamo per le
   lotterie delle feste di paese. Compriamo volentieri il biglietto
  per finanziare l'evento,
  e se poi ci capita di vincere un prosciutto o una bici,
  tanto meglio.\\
  Comunque, tornando a Napoli, non mi sembra che il divieto durò a lungo.
\item[G.] Il problema sorge quando uno gioca al sopra della proprie possibilità,
  o si accanisce con `sistemi' che non hanno senso, o inseguendo numeri 
ritardatari
  o numeri caldi. Perché giocando, mediamente si può solo perdere
  e, come sanno bene gli allibratori clandestini, l'unico modo
  per vincere veramente al lotto è quello di organizzarlo.
\item[A.] E quelli che possono gestirlo legalmente inducono
  la gente a giocare in vari modi. Ora si sono inventati
  il display che mostra il totale vinto a partire da non so quando.
  Io, affianco al totale delle vincite,
  farei mettere per legge anche il totale di quanto
  la gente ha giocato!  
  Ma sull'argomento, oltre all'ignoranza, ci devono essere sicuramente idee errate su quello che veramente dicono le leggi della probabilità, mi sembra di capire. Non credo che i matematici che se ne occupano siano talmente scemi. Ecco il motivo per cui ci eravamo proposte da tempo di approfondire l'argomento. E non solo per mettere in guardia i ragazzi su giochi di questo tipo.
\item[N.] Insegnare queste cose, se non troppo complicate dal punto di vista matematico, ci sembrava più importante di tanti altri dettagli  --  a me almeno così sembrano, o comunque non di fondamentale importanza --  che compaiono nei programmi delle scuole medie, soprattutto alle superiori.
\item[G.] Toh, ma lo sai che la vostra riflessione non è affatto nuova? Qualcuno, un grande del pensiero occidentale, si lamentava della cosa più due secoli e mezzo fa. Anzi, si associava alle lamentele di un altro grande.
\item[A.] Chi sarebbero?
\item[G.] Nientepopodimeno che Hume e Leibniz. Un momento che vi ritrovo quello che scriveva Hume. È una citazione che uso spesso nei miei seminari. Eccola direttamente in inglese, visto che Hume sapeva veramente scrivere bene. E leggetela direttamente voi, altrimenti,
  se ve la leggo io, ``le bon David''
  si rivolterà nella tomba, dopo averlo già fatto per la
  traduzione\footnote{\,A tradurre, ad esempio, ``philosophical speculations'' con ``ipotesi filosifiche''\cite{Hume} ci vuole del corag\-gio, e mi sembra 
    quasi un insulto verso chi 
    voleva 
    fare nella filosofia quello che Newton aveva fatto nella Fisica
    (insomma fondare una sorta di Scienza della Filosofia).}\ldots 
  \begin{quote}
   {\sl  ``The celebrated Monsieur Leibniz has observed it to be a defect in the common systems of logic, that they are very copious when they explain the operations of the understanding in the forming of demonstrations, but are too concise when they treat of probabilities, and those other measures of evidence on which life and action entirely depend, and which are our guides even in most of our philosophical speculations.''}\,\cite{Hume}
    \end{quote}
\item[A.] In effetti\ldots\ se pensiamo a quanto spesso diciamo ``probabilmente questo'', ``questo mi sembra poco probabile'', e così via --
  e non sto parlando solo dei giochi. 
\item[N.] Interessante! Quindi questa lamentela risale a circa tre secoli fa, facendola risalire a Leibniz, e non mi sembra che siano stati fatti sostanziali progressi, almeno nella scuola dell'obbligo. E quindi, dicevo, ci eravamo riproposte di fare qualcosa.
\item[A.] Ma come sai, in giro c'è tanta roba e non sapevamo da dove cominciare. Quindi un prodotto `firmato' ISTAT e appositamente rivolto a bambini e ragazzi ci è sembrato un buon punto di partenza per un lavoro interdisciplinare.
\item[N.] Sì, Alba si sarebbe curata della parte linguistica e letteraria, io di quella matema\-tica, seppur debba confessare che di statistica non ne sappia poi molto.
\item[A.] Va bene, ma quanto vuoi che ne serva? Se il libro è indirizzato addirittura a bambini, un insegnante di matematica se la dovrebbe cavare senza troppe difficoltà, altrimenti c'è decisamente qualcosa che non va, più nel libro che nell'insegnante, direi. 
\item[N.] In effetti la maggior parte delle cose non erano complicate, e anche simpatiche, direi, come quella della svalutazione della moneta dimenticata in una tasca di un cappotto.
\item[A.] Anche se poi secondo me esagerano con la velocità di svalutazione e la storia rischia di generare ansia nel piccolo lettore che ha riposto le sue monetine nel salvadanaio e che dovrà aspettare mesi prima di poterlo aprire. 
\item[N.] Invece quella che ha mandato me in ansia è la fiaba delle streghe, quella che dà il titolo al libro, {\em Le streghe di Bayes}.
\item[A.] Addirittura ansia? Io non ci ho capito niente, ma mi sono rassegnata. Sarà a causa della mia ignoranza, mi son detta.
\item[N.] Sì, forse ansia non è la parola giusta. Il problema è che non ho capito bene cosa volesse insegnare. E poi viene illustrata una strategia decisionale, quella sul tipo di cibo da dare alle streghe in base al colore del cappello, che non mi torna proprio per niente.
\item[G.] Ma nell'appendice c'erano delle indicazioni. Non ti sono servite? 
\item[N.] Peggio! C'è scritto, ecco che lo leggo, ``I bambini \ldots\ verranno guidati alla scoperta del ragionamento bayesiano, che aiuta a cambiare le proprie decisioni quando vengono acquisite nuove conoscenze.''
\item[A.] Mi sembra la scoperta dell'acqua calda. Che non lo sanno tutti che se devi scegliere fra più possibilità la decisione cambia a seconda delle informazioni che hai, ad esempio quanto ti costa fare una cosa piuttosto che un'altra? Se devi andare a fare la spesa in un supermercato e un'amica ti telefona per informarti che in un altro ci sono buone offerte su quello che ti interessa, vai nel secondo, ammesso che non ti costi di più di benzina. Oppure, un classico, se una mattina sei incerta se prendere l'ombrello e senti alla radio che quasi sicuramente pioverà, l'ombrello è meglio che te lo porti! 
\item[N.] Ma puoi portartelo inutilmente, se poi non piove\ldots
\item[A.] È per quello che ho precisato ``quasi sicuramente''. Ma prendiamo un esempio ancora più semplice, senza le incertezze delle previsioni meteo. Se tuo figlio di tre anni sa che metti i biscotti in un certo sportello in basso, si abitua a prenderli lì. Ma se un giorno ti vede
  tornare a casa con
  quelli di cui lui è più goloso e li hai riposti in un cassetto, saprà sicuramente dove
  andarli a cercare, giusto?
\item[N.] Se non è scemo\ldots\ E se non ci arriva, trascina una sedia davanti al cassetto e a modo suo ci si arrampica sopra.
\item[A.] E lo stesso vale anche per gli animali. I gatti di mia sorella, se attraverso l'udito acquisiscono la nuova conoscenza, per usare i termini del libro, che lei sta aprendo una busta di croccantini, si precipitano a raggiungerla.
\item[N.] E già! Altrimenti avrebbero deciso, sempre per usare quei termini altisonanti
  e sebbene `decidere' possa suonare esagerato riferito a gatti,
  di rimanere a poltrire dove stavano. Dico `decidere' perché anche
decidere di non fare niente è pur sempre una decisione. 
\item[G.] A proposito, sentite questa. Una volta mio padre, con la pazienza che aveva, si era divertito a insegnare ai pesci rossi che avevamo in giardino quando stava arrivando con il mangime. Batteva un po' le mani, stando ancora lontano dalla vasca. Una volta mi chiamò per farmi vedere. Era impressionante vedere i pesci che si radu\-navano vicino al bordo
  della fontana da dove avevano imparato che sarebbe comparso!
\item[A.] Certo! Si tratta del ben noto riflesso pavloviano, senza che gatti e pesci debbano conoscere questo Bayes, paese o ``matematico britannico'' che sia, e fare un ``ragionamento bayesiano''. Già, perché questa è un'altra cosa che mi ha confuso: dare il nome di un personaggio che sembra aver fatto cose importanti  --  questo misterioso ``ragionamento bayesiano''  --  a un paese i cui abitanti sembrano a dir poco tonti!
\item[G.] È questa l'impressione che ti sei fatta?
\item[A.] Scusa, questi bayesiani, intesi come gli abitanti di questo ipotetico paesetto di montagna, confidano in quella specie di stregone capo-villaggio e affidano le loro decisioni ai responsi del dio Bias manifestati attraverso l'esito del lancio della pietra. E il loro dio nemmeno ci prende più di tanto. 
E poi, detto per inciso, il nome del paese e quello del dio sembrano fatti apposta per confondere.
\item[G.] Non hai tutti i torti, e ho incontrato pure dei colleghi italiani che confondevano Bayes con il più comune {\em bias} del linguaggio scientifico, che indica una `distorsione sistematica', nel nostro terribile gergo con il quale italianizziamo parole inglesi, specie perché sono più corte e dal significato specifico. Così diciamo che una misura è `baiassata', come parliamo pure di
  esperimenti a `targhetta fissa',
  orribile traduzione di {\em fixed target},
  che sarebbe `bersaglio fisso'..
\item[A.] Mamma mia!
\item[G.] Quindi, tornando alla fiaba, i bayesiani ti sembra che
  vengano presentati nella storia come dei sempliciotti?
\item[N.] Sicuramente poco intelligenti. Alla fine si viene pure a scoprire, mediante l'arcano della filastrocca, che alle streghe con il cappello nero piacciono solo cibi salati. Ci voleva tanto a capirlo, visto che ogni volta che mettevano cibi dolci nel cappello nero succedeva il pandemonio? E nella fiaba c'è scritto che la storia andava avanti da un bel po' e che i poveri bayesiani avevano provato diverse soluzioni. Come minimo non ripetere quello che ti ha dato un risultato indesiderato e continua con quello con cui ti è andata bene!
\item[G.] Sì, come la barzelletta di quello al distributore automatico che seguitava a inserire monete e raccogliere merendine, dicendo a quelli che erano in coda ``scusate, ma finché vinco continuo a giocare io!''
\item[A.] O tipo il tacchino induttivista di Russel\ldots
\item[G.] Già, più o meno.
\item[N.] Aspe', me lo ricordate, questo tacchino filosofico?
\item[A.] Ogni mattina, quando sentiva i passi della padrona, il tacchino rafforzava il proprio convincimento che la signora stava arrivando per portargli da
  mangiare\ldots\ perché gli voleva bene, pensava lui\ldots
\item[G.] Finché il giorno prima della festa la padrona venne per tirargli il collo.
\item[A.] E poi, scusate, al di là degli aspetti matematici, che magari vi chiarirete voi due, e tornando alla ``filastrocca numerica'', come viene definita nell'appendice\ldots
\item[G.] Filastrocca numerica? Veramente la chiamano così? Me ne ero dimenticato, o questo dettaglio mi era sfuggito.
\item[A.] Una filastrocca numerica, dal punto di linguistico,
  è una filastrocca fatta di numeri, non so\ldots\ non me ne viene nessuna\ldots
\item[N.] Forse quella della tabellina del sette, che gira dalle nostre parti,
\begin{quote}
 {\sl Sette, quattordici, ventuno, ventotto, \\
è cascata a moje 'e bassotto\\
 e ha fatto un grande botto. \\
Sette, quattordici, ventuno, ventotto.}
\end{quote}
\item[G.] Ah, ricordo! Me la ripeteva mia nonna. E funzionava talmente bene che alle elementari l'inizio della tabellina del sette lo sapevano tutti. I guai venivano da trentacinque in poi\ldots  \ 
Ma scusa, Alba, le tue obiezioni erano solo se la filastrocca della fiaba si può definire numerica?
\item[A.] Certamente no, anzi nemmeno mi ero posta il problema, concentrata com'ero a capire la sostanza. Mi riferivo al fatto che, come c'è scritto nell'appendice, grazie a tale filastrocca ``Nora riuscirà a fare delle previsioni più affidabili''. A me sembra che Nora non faccia nessun tipo di previsione, ma che cerchi di indo\-vinare il gusto della strega basandosi sul colore del cappello, e non sono sicura che si tratti della stessa cosa. Previsione sarebbe stato tentare di indovinare il colore del cappello del giorno dopo.
\item[N.] A proposito, nel testo a un certo punto si parla di una previsione del genere. Aspettate che lo ritrovo. Ecco, in fondo a pagina 49, e si riferisce a Nora: ``si incuriosiva nel pensare quale cappello la mattina avrebbero posato fuori dall'antro: nero o viola?''
\item[A.] E questa sì che sarebbe stata una previsione! Ma poi la cosa muore lì, mentre pensavo che sarebbe rispuntata da qualche parte nella fiaba. 
\item[G.] In effetti quella sarebbe stata veramente
  una previsione probabilistica interessante,
  insomma non banale, rispetto a quella di cercar di capire il gusto dal
  colore del cappello.
\item[N.] Dici banale in quanto, con le ipotesi della fiaba, tale
  probabilità è costante
 e non dipende da quanto è successo nei giorni precedenti,  giusto?
\item[G.] Eh sì, è come nel lancio di un dato regolare.
  Se volete, anche dire che ciascuna faccia ha probabilità 1/6
  è una previsione probabilistica. Ma è banale in quanto dai lanci precedenti
  non si impara nulla sui successivi. Ben diverso sarebbe il caso
  di un solido esagonale irregolare, tutto storto e per giunta fatto di un
  materiale non omogeneo.
  Nel caso delle streghe, 
  prevedere, seppur probabilisticamente, il colore del cappello
  del mattino dopo sarebbe stato -- quello sì! -- 
  un bello spunto per un problema veramente inferenziale-predittivo. 
\item[A.] Inferenziale-predittivo? Ora ti ci metti pure tu. Mi ricorda una battuta di {\em Smetto quando voglio}, di quello che si dichiarava esperto in statistica
  inferenziale.\footnote{\,\url{https://www.youtube.com/watch?v=seEhOShK0cc} (trailer, da 0:35).}
  Ho rivisto il primo episodio poco tempo fa, ma non ho nessuna idea di cosa sia questa statistica. Una mia amica, che durante il film, qui a casa, chissà a cosa stava pensando, aveva addirittura capito `statistica influenzale' e diceva che forse quel tipo faceva le statistiche di quelli a letto con l'influenza. 
\item[G.] Geniale la statistica influenzale!
\item[N.] Ah ah, sì, divertente\ldots\ Ma tornando alla fiaba, scusa Giulio, quale sarebbe il problema inferenziale-predittivo che si potrebbe inserire nella storia? Io, come dicevo, di probabilità e statistica so veramente i rudimenti: casi favorevoli e casi possibili, media, deviazione standard, binomiale, gaussiana, \ldots 
\item[A.] \ldots\ che per me già è arabo. Io l'unica cosa che ho imparato, dalle elementari all'università, che avesse a che fare con la statistica è la famosa poesia di Trilussa.
\item[N.] E io poco più\ldots\ veramente l'ABC, quelle poche cose che ci hanno fatto a biologia. Ci potresti fare un esempio, possibilmente legato alla fiaba delle streghe? Intendo un esempio di problema inferenziale-predittivo.
\item[G.] Era quello che avevo in mente: inferire --  insomma `cercare di indovinare al meglio' --  quante sono le streghe dei due gruppi, cappelli neri e cappelli viola; cercare di prevedere il colore del cappello che sarebbe stato esposto all'ingresso della caverna una certa mattina.
\item[A.] Beh, se dobbiamo indovinare, ci possiamo inventare quello che ci piace, anche se già vedo un problema. Verificare se abbiamo previsto correttamente il colore del cappello è facile: basta andare il giorno dopo davanti alla grotta e lo vediamo  --  scusate, mi sto mettendo nella storia\ldots\ 
Ma verificare il numero delle streghe con il cappello nero e di quelle con il cappello viola potrebbe non essere possibile. Non credo ci farebbero entrare e permetterci di contarle, tutte in fila, ciascuna col suo bel cappello. Forse potremmo chieder loro. Ma immagina se ci rispondono, o se ci dicono i numeri veri, dispettose come sono!
\item[N.] O addirittura se ne potrebbero andare da quella caverna, di notte, per recarsi altrove, e non lo sapremo mai\ldots
\item[G.] Osservazioni interessanti. Preciso che per `indovinare' intendevo esprimerci in termini probabilistici. E non solo qualitativamente, come potrebbe essere ``è più probabile viola che nero'', ma anche in termini quantitativi, ad esempio ``viola è cinque volte più probabile di nero'', e così via.
\item[A.] Ma in entrambi i casi, sia riguardo il numero di cappelli viola che riguardo il cappello che verrà esposto fuori dalla caverna, stiamo cercando di `indovinare', seppur in termini probabilistici, come dici tu. Quindi in cosa differisce l'aspetto inferenziale da quello predittivo?
\item[G.] Sì, è una sottigliezza, se vogliamo, perché in entrambi i casi si tratta di quantificare quanto una possibilità è più probabile di un'altra. Ma, se ci pensi, una piccola differenza che può essere filosoficamente importante  --  scusate il parolone  --  l'hai già fatta notare tu. Il colore del cappello esposto davanti alla caverna è un `fatto' osservabile. Basta aspettare e lo sapremo. Ma per qualche motivo potremmo non accontentarci semplicemente di aspettare. Potremmo invece essere interessati a anticipare se è più facile che si verifichi nero o viola, ben consci che si tratta di un nostro giudizio probabilistico. Fa parte dell'animo umano. Come disse qualcuno, ``meglio fare previsioni, seppur con incertezza, che non predire per
  niente''.\footnote{\,``Mieux vaut
    prévoir sans certitude que de ne pas prévoir du
    tout''\,\cite{ScienceHypothese}
  (``It
  is far better to predict without certainty, than never
to have predicted at all'').} 
\item[N.] Interessante, chi era?
\item[G.] Credo Poincaré\ldots\ Henry Poincaré.
\item[A.] ``Chi era costui?'', mi verrebbe da dire\ldots\ Scusate, continuate\ldots
\item[N.] Sì, prevedere con incertezza, come piove o non piove domenica prossima. Anche se non siamo sicuri al 100\%, ci è senz'altro utile avere un'idea se è più facile che piova piuttosto che sia sereno. Almeno possiamo ridurre il rischio legato a eventi spiacevoli. Poi, magari, se siamo particolarmente interessati all'evento, seguiremo le previsioni durante la settimana, perché a mano a mano che ci avviciniamo al giorno di interesse saranno più accurate -- a parte mettersi d'accordo su cosa significhi `pioggia', insomma quanti millimetri di acqua ci servono per definirla tale, ma lasciamo stare queste sottigliezze.
\item[A.] Sì, come quelli che disquisivano su quale fosse il granello di sabbia che trasforma un `non mucchio' in `mucchio'. Diciamo che ci siamo capiti. Non perdiamoci in sofismi, anche se capisco bene che in certi casi un po' di precisione ci vuole.
\item[G.] Certo, un grado di precisione adeguato allo scopo. E comunque potremmo sempre ampliare lo spazio delle possibilità: meno di un millimetro, da due a cinque millimetri, e così via, precisando che ci stiamo riferendo a una ipotetica provetta posta in mezzo alla piazza del paese, se proprio vogliamo essere pignoli. Ma lasciamo stare queste sottigliezze e concentriamoci sul semplice caso della fiaba: nero o viola, e senza nessun effetto ottico che ci possa far confondere i due colori.
\item[A.] Ok, per quello che riguarda il problema predittivo mi sembra che ci siamo capiti. E per quanto riguarda quello inferenziale?
\item[G.] Nel caso della fiaba, o almeno della variante che stiamo prendendo in considera\-zione, è quando cerchiamo di capire la proporzione di streghe dal cappello viola, in quanto non è qualcosa che possiamo osservare direttamente, almeno in linea di principio  --  se poi ciò fosse possibile, magari in seguito, tanto meglio. 
\item[N.] Mi sembra di aver capito. Possiamo solo osservare gli effetti di una proporzione più o meno alta di tali streghe. Più essa è alta più facilmente ci aspettiamo di osservare un cappello viola fuori dalla grotta.
  E, ragionando al contrario, dall'osservazione di tanti cappelli viola e di pochi neri
  possiamo dedurre che le streghe dal cappello viola siano la maggioranza\ldots
\item[G.] Anche se non si tratta propriamente di deduzione.
\item[A.] Come no? Non è come quando un detective, analizzando
  la scena del delitto, deduce chi può essere stato a compiere il crimine?
  Anche qui abbiamo una relazione cause-effetti e dagli effetti cerchiamo di risalire
  alle cause.
\item[G.] Sì, il concetto è quello. È il nome del processo logico che non è corretto.
  Non si tratta di deduzione, anche se è quello il termine che usa l'inventore
  del più famoso investigatore.
\item[A.] Intendi Sherlock Holmes?
\item[G.] Per l'appunto. Non ti ho mai parlato di 
  {\em Plato and a platypus walk into a bar}?\,\cite{PlatoPlatypus}
  Mi sembra impossibile! Quando l'ho scoperto ne ero rimasto talmente
  entusiasta che ci ho addirittura scritto una recensione per 
  {\em Scienza per Tutti}\,\footnote{\,\url{http://scienzapertutti.infn.it/rubriche/un-libro-al-mese/}, settembre 2014.}
  dei Laboratori Nazionali di Frascati. 
  E ci sono proprio alcune pagine divertenti che chiariscono
  il concetto  mediante
  una efficace barzelletta, proprio su  Sherlock Holmes
  e l'immancabile Dr. Watson\footnote{\,\url{http://www.roma1.infn.it/~dagos/PlatoPlatypus_30-31.pdf}.} -- l'intero libro è pieno di barzellette
  che spiegano la filosofia, e lo fanno spesso meglio di manuali scolastici.
  Poi te lo presto. 
\item[N.] Quindi, cerco di ridire correttamente quello che avevo in mente,
  usando il termine `inferenza'\ldots \ 
   Dall'osservazione di tanti cappelli viola e di pochi neri
   possiamo inferire che le streghe dal cappello viola siano la maggioranza.
   Ma farlo quantitativamente, ovvero esprimermi in termini probabilistici
   sulle possibili composizioni di cappelli neri e viola, non ne sarei in grado.
 \item[A.] Io pure non riuscirei ad andare oltre l'affermazione qualitativa
   di Noemi. 
   Ma mi verrebbe da dire che, se un giorno osservo un cappello nero e il giorno dopo un cappello viola, sto al punto di partenza.
\item[G.] Perfetto! È inoltre interessante l'idea di Noemi di chiamare `effetto' il colore del cappello osservato. Difatti problemi di questo tipo vengono spesso modellizzati in termini di cause e effetti, anche se a rigore, nel caso generale,
  non si tratta di causalità, bensì di condizionamento.
\item[N.] Un po' come quando si dice che ``{\em correlation does not imply causation}''?
\item[A.] Sarebbe?
\item[N.] Immagina gli incidenti di bagnanti in mare e le vendite di gelati. Nelle settimane in cui aumenta il numero di gelati venduti aumenta anche il numero di incidenti in acqua, ma non perché essi siano causati dalle vendite di gelati, a parte quei casi di chi si butta in acqua dopo essersi mangiato due chili di cassata.
\item[A.] Ho capito. È la settimana calda a essere `causa' di entrambi, anche se essa non è la sola causa, come si può ben capire.
\item[G.] Sì, grazie, era quello che intendevo, e senza parlare
  delle correlazioni spurie.
\item[N.] Sì, c'è addirittura qualcuno che si è
  divertito a raccoglierne un certo numero \cite{SpuriousCorrelations},
  come quella dell'importazione di greggio degli Stati Uniti
  dalla Norvegia e le morti di automobilisti in incidenti con
  treni.\footnote{\,\url{http://www.tylervigen.com/spurious-correlations}\,.}
\end{description}
  \begin{center}
\epsfig{file=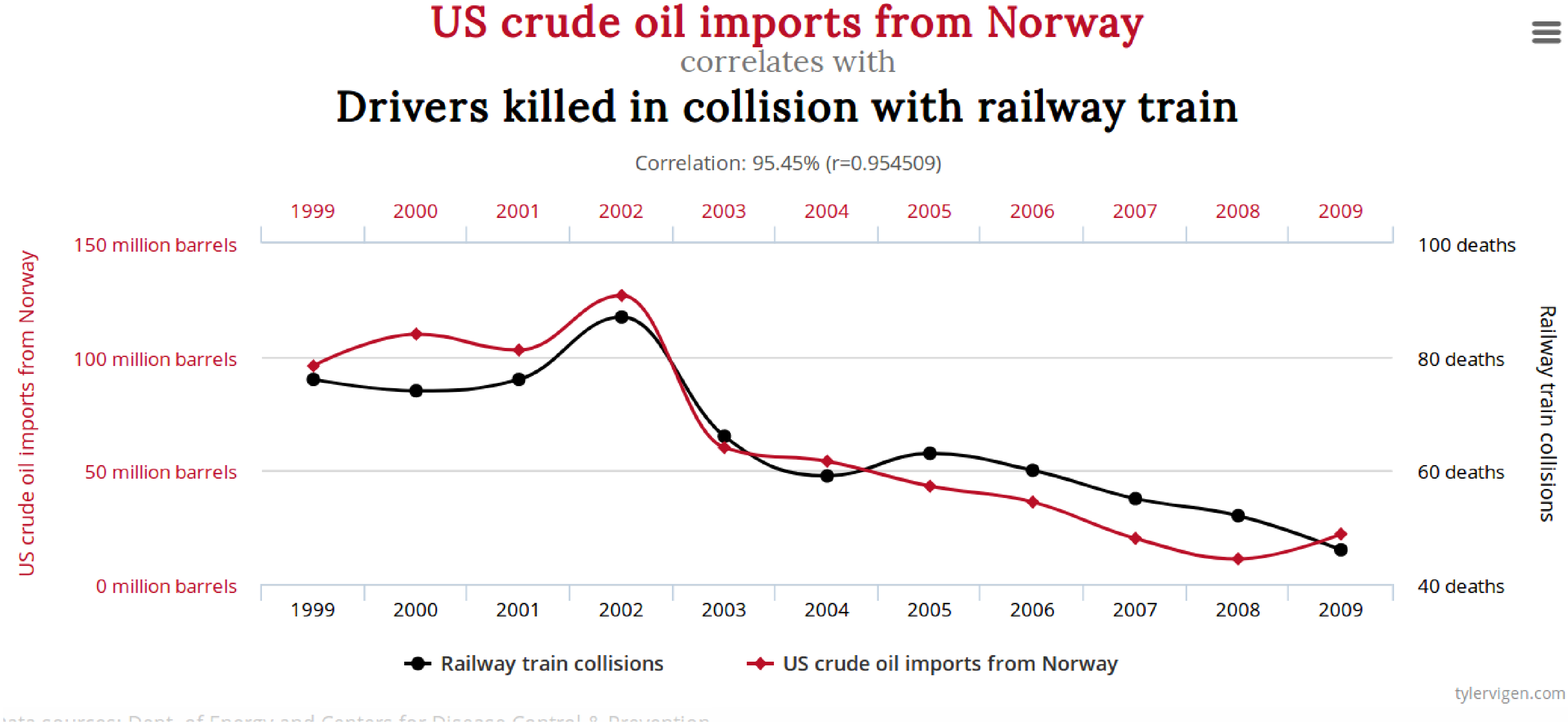,clip=,width=\linewidth}
  \end{center}
\begin{description}
\setlength\itemsep{0.5mm} 
\item[G.] Sì, non solo divertenti (mi riferisco alle correlazioni, non agli incidenti di per sé), ma anche istruttive, soprattutto se pensiamo a
 quando sentiamo dai media che una cosa ci fa male e dopo un paio di anni, se non mesi, ci fa addirittura bene. Vatti a fidare!
\item[N.]  In effetti non è facile trovare nessi causali degli agenti chimici e biologici
  su una macchina complicata come il corpo umano.
\item[G.] D'accordo, ma dovrebbero essere come minimo più prudenti
  prima di fare certe affermazioni, soprattutto perché basate su certi metodi statistici
  a dir poco fantasiosi.\footnote{\,Chi desidera approfondire veda i par. 2 e 3 di \cite{Waves-sigmas} (pp. 3-10) e le citazioni ivi contenute.}
  Ma per adesso continuiamo con lo schema cause-effetti. Un effetto è quello che possiamo vedere direttamente, come il colore del cappello, che discende, seppur in modo non deterministico, da una `causa', in questo caso la proporzione di cappelli viola.
\item[N.] Scusa, ma a rigore la causa sarebbe la volontà delle streghe a esporre il cappello di una di loro scelta a caso, ma capisco cosa intendi dire.
\item[G.] Sì, l'osservazione è pertinente, grazie. In effetti abbiamo delle concause, come in tutti gli eventi che ci accadono nella vita.
  Era quello che osservava Alba un momento fa: anche nel caso della correlazione
  fra vendite di gelati e incidenti in mare la causa comune sarà
  pure l'estate, ma se uno al mare non ci va il rischio di affogare
  non lo corre di certo. E, tornando a noi, 
  qual è la `causa' per la quale oggi siamo qui, noi tre,
  a casa sua a discutere di queste cose? Veramente verrebbe da dire `infinite'. 
Tornando alla fiaba, è chiaro che il cappello appare fuori dell'ingresso della grotta perché le streghe vogliono così. Ma l'effetto nero o viola è legato, seppur probabilisticamente, alla proporzione di streghe dal cappello viola. Quindi l'altra causa, che addirittura potremmo chiamare `primaria', la possiamo ignorare, ai fini della nostra inferenza, in quanto la diamo per scontata.
\item[A.] Quindi  --  fatemelo dire a modo mio per vedere se ho capito  --  stiamo dicendo che lo stesso effetto, ad esempio `cappello nero', può discendere da tante `cause' possibili, che sono tutte le proporzioni di streghe dai cappelli viola che possiamo prendere in considerazione.
\item[G.] A questo punto farei un diagramma per rappresentare
  il modello logico che stiamo prendendo in considerazione. Mi vien da alzarmi
  e dirigermi verso una lavagna che qui manca...\  -- aberrazione professionale, scusate.
\item[A.] Lavagne a casa non ne ho di certo, mi bastano 
quelle della scuola! Ma da scrivere ne ho quanto ne vuoi. Un attimo che
  provvedo.
\end{description}
%
\begin{center}
{\Large  $[\,$ Pausa $\,]$ }
\end{center}  
\mbox{}\\
\begin{description}
  \setlength\itemsep{0.5mm}
\item[A.] Eccomi, con fogli, penne e matite. E ho portato
  anche dell'acqua\ldots\ da bere, non perché penso che Giulio
 abbia voglia di mettersi a  fare acquarelli\ldots 
\item[G.] Potrebbe essere un'idea\ldots\  ma non riuscirei mai a eguagliare
i bei disegni del libricino dell'ISTAT -- veramente ben fatti.
\item[A.] Ho messo pure a bollire acqua per farci un the o una tisana.
  Volete che prepari pure un caffè?
\item[G.] No, grazie, acqua e una qualche tisana mi vanno bene a quest'ora.
\item[N.] A me basta l'acqua, grazie.
\item[G.] Allora, ripartiamo da quello che stava dicendo Alba, 
  che lo stesso effetto, ad esempio `cappello nero'
  può discendere da tante `cause' possibili, che sono tutti
i possibili numeri di streghe dai cappelli viola
rispetto al totale.
  E lo stesso vale per `cappello viola'. 
  Questo è il diagramma
  del nostro modello, limitatamente a un solo giorno e nel quale
  mettiamo delle frazioni inventate, ma utili per
  afferrare il tipo di ragionamento. 
  \begin{center}
  \epsfig{file=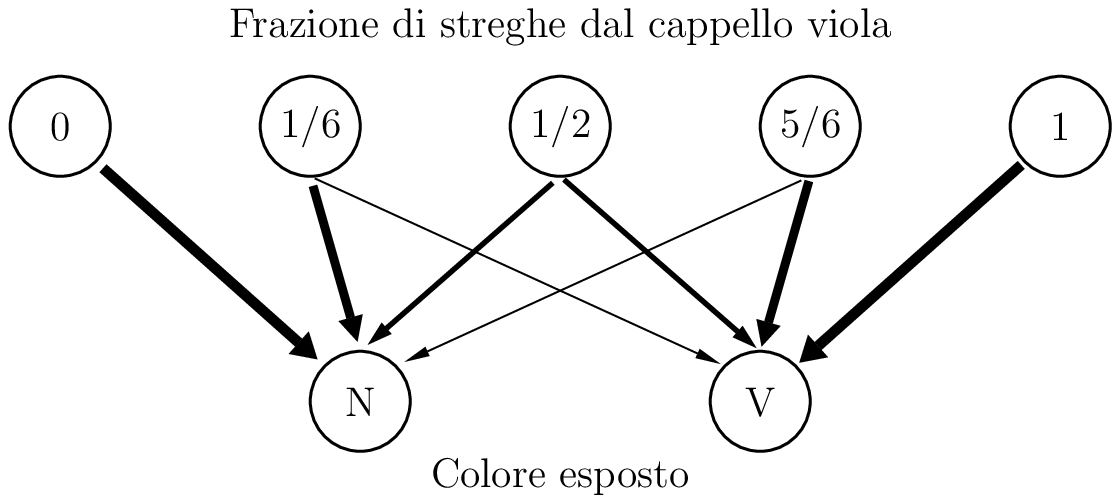,clip=,width=0.7\linewidth}
  \end{center}
  Nel diagramma lo spessore delle frecce dà un'idea
dei diversi valori probabilità.
\item[A.] Come mai hai messo 1/6 e 5/6, che non compaiono nella fiaba?
\item[G.] Dicevo che sono comunque valori fittizi. Li ho scelti
 perché sono tipici di probabilità nel lancio di dadi,
  con i quali abbiamo tutti una certa familiarità.
  Quindi, ad esempio, la frazione `$1/6$' dà Nero esattamente
  con la stessa probabilità con la quale un dado può mostrare
  una particolare faccia.
\item[N.] E quella indicata con 1/2 corrisponde, ad esempio, all'insieme
delle facce con numeri pari, e così via.
\item[A.] Questo diagramma somiglia alle mappe concettuali che
  usiamo a scuola in diverse materie, per connettere vari concetti,
  argomenti o avvenimenti.
\item[N.] Ma in questo caso i collegamenti sono probabilistici, 
  anche se alcuni casi possono essere deterministici, come quelli 
  a partire da `0` e da `1'. Ma questi vanno visti come casi 
 speciali.
\item[G.] E, dato il loro uso nel cosiddetto ragionamento bayesiano,
  grafici di questo tipo vengono chiamati
  reti bayesiane, {\em Bayesian networks} in inglese.
\item[N.] Quindi il problema inferenziale di uno schema del genere 
  sarebbe
  quello di valu\-tare la probabilità delle diverse frazioni
  a partire dal cappello osservato. Se osservo Nero,
  è più facile che provenga da `0', e sicuramente
  non da `1', che corrisponde a soli cappelli viola. 
\item[G.] Sì, ma solo se ritenevi le quattro 
  cause che possono aver dato Nero
ugualmente possibili. Se per esempio all'ipotesi '0'
  non ci credi, perché ad esempio hai già osservato 
  un cappello viola, seguiti a non crederci, anche se è quella che
  produce solo cappelli neri. Così pure, se per qualche motivo
  credevi a `5/6' molto di più che a `1/6',  
  seguiti a crederci di più, anche se un po' meno di prima.
\item[N.] Quindi a mano a mano che osserviamo cappelli neri
  le nostre credenze sulle varie ipotesi si evolvono, per così dire.
\item[G.] Ben detto. Bisogna solo capire come.
\item[A.] Ma a me questa storia delle `credenze` non è che mi piaccia
  molto! Se adesso dobbiamo dar retta a tutte le credenze che si sentono in
  giro abbiamo finito!
\item[N.] Ma si tratta di credenze razionali, basate su fatti
  e su ragionamenti. Non  dogmi, fantasie o allucinazioni.
\item[G.] So bene che molti diffidano della parola `credenza',
  in quanto le attribuiscono connotati irrazionali o non scientifici.
  L'importante è capire di cosa stiamo parlando, anche perché è inevitabile
  che il linguaggio umano abbia ambiguità e che le stesse parole abbiano
  molteplici significati.
\item[A.] Certo, e di sicuro non avevamo minimamente pensato alla credenza
  della cucina\ldots 
\item[G.] E la probabilità sta appunto a quantificare il grado di credenza 
  razionale su qualche cosa,
  quello che in inglese viene detto papale papale {\em degree of belief}.
  Ma sembra che in Italia si abbia particolarmente paura 
della traduzione letterale di {\em belief}, e ``{\em degree of belief\,}'' 
viene comunemente reso come ``grado di fiducia''.
\item[N.] E quindi, riassumendo, il problema inferenziale è, in questo caso, quello di effettuare affermazioni probabilistiche su ciascuna proporzione possibile a partire dagli effetti osservati, ovvero dai colori dei cappelli che abbiamo visto, o che ci sono stati riferiti, giorno dopo giorno.  
\item[G.] Sì, e il diagramma di sopra si modifica in quanto 
          non abbiamo più soltanto N e V, ma una particolare 
          sequenza di colori, ad esempio $NNVNN\ldots$, e così via.
          Quindi il diagramma diventa
           \begin{center}
           \epsfig{file=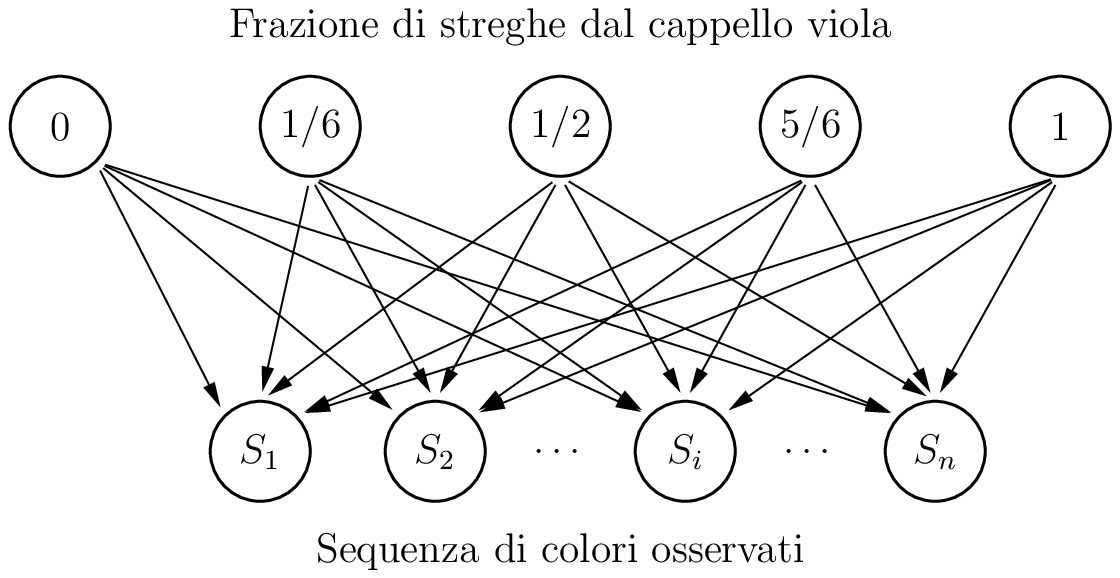,clip=,width=0.7\linewidth}
           \end{center}
          ove con $S_i$ ho indicato la generica sequenza (`$i$-esima', come si dice). 
\item[N.] E dalla particolare sequenza osservata possiamo 
          valutare quale frazione di cappelli viola è più o meno credibile.
\item[A.] Ecco, ora comincia a avere un senso la frase dell'appendice ``che aiuta a cambiare le proprie decisioni quando vengono acquisite nuove conoscenze'', riferita al ``ragionamento bayesiano''.
\item[G.] Uhm, eppure quella frase contiene a mio giudizio un bel pasticcio, perché mischia valutazioni di probabilità con problemi decisionali. In effetti, il ragionamento bayesiano aiuta, a rigore,
  solo a cambiare le probabilità delle varie
  ipotesi alla luce delle nuove informazioni.
Una decisione è una cosa più complicata perché dipende sia dalle probabilità dei possibili effetti che dalla loro deside\-rabilità, ovvero dal `valore',
spesso assolutamente soggettivo, che assegniamo a ciascuno di essi. Ecco perché se le previsioni meteo danno un bel temporale, seppur con bassa probabilità, adottia\-mo tutte le precauzioni necessarie: non soltanto prendere l'ombrello, ma anche indossare scarpe adeguate, eventualmente prendere la macchina, o addirittura rimandare l'uscita, se non strettamente necessaria. 
\item[N.] Mi sembra di aver capito, ma non saprei da dove cominciare. La probabilità la so calcolare solo in alcuni casi facili.
\item[G.] Sì, lo immagino. E qui non si tratta di valutarla da zero, per così dire, ma di aggiornarla.
\item[A.] Aspettate, io che ci capisco poco ora mi sono persa completamente. Ci stai dicendo che per calcolare una probabilità dobbiamo partire da un valore di probabilità precedente che dobbiamo aver già calcolato, giusto?
\item[G.] Giusto.
\item[A.] E magari per calcolare quello precedente ce ne serve uno ancora precedente\ldots\  Ma allora la storia non finisce più! Mi ricorda quelli che credevano che la Terra fosse poggiata su una tartaruga, che stava sopra un'altra tartaruga, e giù all'infinito. O qualcosa del genere.
\item[N.] E già, scusa, ma una qualche probabilità di partenza ci deve essere, e quella come la calcoliamo? Mi sembra che le obiezioni di Alba siano più che legittime.
\item[G.] Sì, capisco bene che la storia possa ricordare, per fare un esempio più terra terra, la scena di Fantozzi e il ragionier Filini sulla barca, impigliati nella rete, con Gesù che si avvicina camminando sull'acqua del
  lago:\footnote{\,\url{https://www.youtube.com/watch?v=x0U4eE7fqmY} 8:15 }
{\sl  \begin{itemize} 
 \item[]  Avete pesci?
 \item[]  Non abbiam pescato, dottore.
 \item[]  Avete pane?
 \item[]  No, doveva portarlo lui\ldots 
 \item[]  E allora che me moltiplico io?
 \end{itemize} }
\item[A.] Mitica! Ma ce l'hai ricordata per distrarci o perché c'è una qualche analogia?
\item[G.] Una qualche analogia c'è, se proprio vogliamo. Ma `Bayes', insomma la formula che porta il suo nome, non moltiplica pani e pesci bensì probabilità, anche se non si tratta propriamente di moltiplicazione, in quanto il fattore di moltiplicazione dipende dalla probabilità che si vuole moltiplicare.
\item[A.] Pretendi che capisca? 
\item[G.] Detto così sarebbe un miracolo, per rimanere in tema. Facciamo un passetto alla volta.
\item[N.] Con l'ultima cosa che hai detto mi sono confusa di nuovo pure io perché una rapida ricerca sul teorema di Bayes l'avevo fatta, e a vedere la formula sembrava proprio una normale moltiplicazione.
\item[G.] Il problema è che quel fattore moltiplicativo è scritto sotto forma di una frazione, la quale ha, a denominatore, una probabilità che, esplicitata opportunamente, contiene la probabilità che si vuole aggiornare\ldots
\item[A.] Spero che almeno voi vi stiate capendo, visto che non so nemmeno di quale formula state parlando, se non me la fate vedere, e ammesso che la capisca\ldots 
\item[G.] Sì, Alba, scusa, hai ragione. E il motivo è che pure io non so esattamente in che direzione andare. Non è che io sia venuto qui con una lezione pronta, con tanto di lucidi. Cerco di procedere `a vista', in base a quello che mi dite e a come reagite a quello che vi dico. Ma facciamo così, lasciamo stare per ora la formula standard che si trova in giro  --  eventualmente ci torniamo  dopo --  e partiamo da quella che ci insegna come aggiornare il rapporto delle probabilità di due ipotesi mutuamente esclusive. 
\item[N.] Tanto se abbiamo il nuovo rapporto di probabilità è facile calcolare la nuova probabilità di ciascuna di esse. 
\item[G.] Il vantaggio è che il fattore di aggiornamento del rapporto è veramente un `fattore esterno', perché non dipende dalle probabilità di partenza.
\item[A.] Quindi, l'analogia con la scena tragicomica di Fantozzi potrebbe essere che, invece di pani e pesci, `Bayes', per così dire, moltiplica i rapporti di probabilità.
\item[G.] E quelli ce li devi avere --  almeno uno, diciamo, ovvero devi avere almeno due ipotesi e devi esserti fatta un'idea di quanto una sia più credibile dell'altra. Altrimenti un ipotetico `Bayes' fantozziano direbbe, andandosene, ``e che me moltiplico io?''.
\item[A.] Ok, chiaro. Quindi, ad esempio, io ipotizzo due possibilità e dico che, in base a mie considerazioni che non ti sto a giustificare, una mi sembra dieci volte più probabile dell'altra.
  E ora arriva  --  quasi quasi mi sembra di vederlo, camminando lentamente sull'acqua\ldots\ --
  questo Bayes e -- tac! -- trasforma il mio dieci in trenta, avendolo moltiplicato per tre.
\item[G.] Molto suggestiva l'immagine di questo Bayes fantozziano
 che compie la moltiplicazione dei rapporti di probabilità. Ma diciamo che dà l'idea, salvo aggiungere qualche precisazione\ldots 
\item[N.] Ma se io, per motivi miei, pensavo inizialmente che il rapporto di probabilità valesse due, Bayes me lo farebbe diventare sei, giusto?
\item[G.] Diciamo che è giusto, ma solo se avete le stesse informazioni empiriche e lo stesso modello di come gli effetti dipendono probabilisticamente dalle cause.
  Ecco, nel caso più generale, con tante cause e tanti effetti, abbiamo un
  modello del genere: 
   \begin{center}
           \epsfig{file=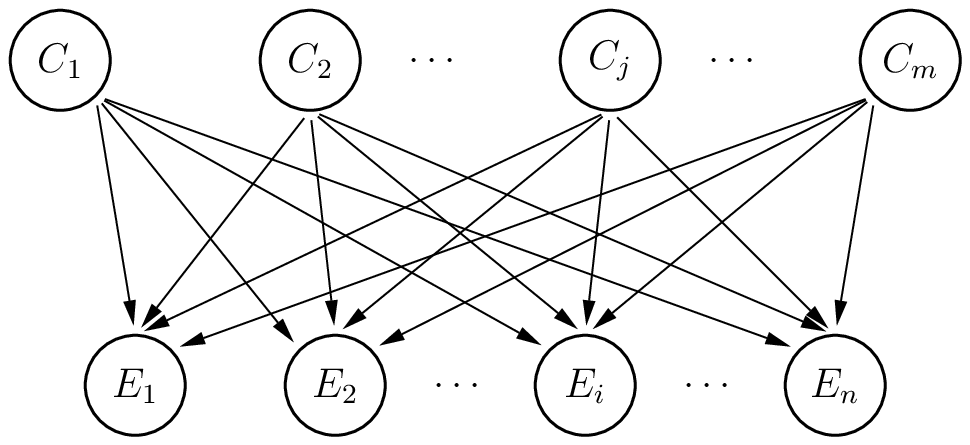,clip=,width=0.7\linewidth}
   \end{center}
Ma torniamo alla fiaba delle streghe e cerchiamo di metterci d'accordo su quali sono le ipotesi di lavoro. Incominciamo quindi prendendo per buone le informazioni contenute nella filastrocca: alle streghe con il cappello nero piace solo salato, mentre alle streghe con il cappello viola?
\item[N.] A sei su sette piace dolce. E forse potremmo disegnare
  un diagramma del genere anche per queste probabilità. Posso? E provo anche a fare lo spessore
delle frecce che diano un'idea delle probabilità. Eccolo
 \begin{center}
           \epsfig{file=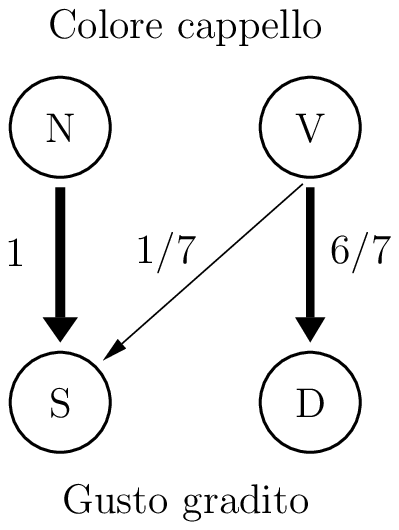,clip=,width=0.28\linewidth}
   \end{center}
\item[G.] Ben fatto! E anche su questa rete bayesiana minimalista 
si potrebbe creare un problemino inferenziale: valutare, 
sapendo il tipo di cibo offerto e che la strega lo ha gradito,
se il cappello era nero o viola. 
\item[A.] Se è stato offerto un dolce e la strega ha gradito, 
si tratta sicuramente di una strega dal cappello viola, 
come è chiaro dal diagramma. Uhm, mica una brutta idea queste
reti bayesiane per visualizzare questo tipo di situazioni!
E forse pure altre\ldots\ sto pensando ai test di inferenza delle prove Invalsi,
sui quali i ragazzi si dimostrano particolarmente deboli.
Chissà se queste rappresentazioni grafiche li possono aiutare, 
ci devo riflettere. Scusate stavo pensando a un problema scolastico\ldots\ uno
dei tanti\ldots\ e ho perso il filo. Stavamo dicendo?
\item[G.] Dicevi di cosa inferire se veniamo a sapere che alla strega
  fosse andato bene Dolce. 
  E nel caso avesse gradito Salato?
\item[A.] Direi che è più probabile che si tratti di cappello nero.
\item[G.] E tu che ne pensi?
\item[N.] Per quanto detto prima, mi verrebbe da dire che deve dipendere
anche da quanti cappelli neri e viola ci sono nella grotta. 
Se, al limite, avessimo mille streghe, di cui solo una con il cappello nero,
anche in questo caso scommetterei che si tratta di cappello viola.
\item[A.] Stai dicendo che anche se le offro Salato
  e non si arrabbia devo sempre pensare a una strega dal cappello
viola? 
\item[N.] Ma anche se si dovesse arrabbiare! Semplicemente perché al 99.9\%
  verrà esposto un cappello viola. E quindi punterei
su Viola qualsiasi cosa accada.
  Sarebbe simpatico risolvere il problema anche quando non siamo in questo caso limite di streghe Viola a Nero 999 a 1. Ad esempio solo uno a uno o due a uno.
\item[A.] Il numero di streghe nella grotta è ventuno, 
a una su sette di quelle con il cappello
viola piace salato e quindi 
esse possono essere
solo o sette, o quattordici o ventuno\footnote{\,Qui ci potrebbe essere una ambiguità, in quanto 1/7 si potrebbe riferire all’intera ‘popolazione’ di streghe dal cappello viola, delle quali un numero imprecisato si troverebbe nella grotta. In questo caso le possibilità del numero di streghe dal cappello viola non sarebbero soltanto 7 o 14, bensì tutti i numeri da 1 a 20 (escludendo 0 e 21 una volta che sappiamo che ci sono cappelli di entrambi i colori). Per semplicità scegliamo l’interpretazione che 1/7, esatto, si riferisca al numero di streghe dal cappello viola presenti nella caverna.}\ldots\ mi sembra la filastrocca 
dei nostri nonni, ah ah.
\item[N.] Già, sono ventuno. C'è scritto nella fiaba, ma ora che ci ripenso
ma non si capisce perché abbiano dato il numero preciso, 
visto che nella storia non viene usato.
\item[G.] Ma a noi ci fa comodo per fissare le idee e lo usiamo.
\item[N.] Ma ventuno streghe dal cappello viola non può essere, perché se così fosse non verrebbe mai fuori il cappello nero. Mentre tale colore è stato sicuramente osservato.
  Così come non può essere zero.
\item[A.] Quindi, non essendo né zero né ventuno abbiamo solo due possibilità, sette o quattordici, dico bene?
\item[G.] E fra sette e quattordici, su quale possibilità scommettereste? Insomma è più probabile che siano sette o quattordici?  
\item[N.] A priori direi che sono equiprobabili, nel senso che non ho alcun motivo razionale per credere a una possibilità piuttosto che a un'altra. Successivamente,
  come abbiamo detto, le probabilità 
  dipenderanno dalla sequenza Nero/Viola osservata nel passato.
\item[G.] E Viola sembrerebbe essere stato dominante, anche se non viene detto
  esplicitamente.
\item[A.] Come fai a dirlo?
\item[N.] Da qualce parte c'era scritto.   
  Mi passi un attimo il libricino?
  Ecco, in alto a pagina 50 si dice che ``alla bambina sembrava che i cappelli viola fossero di più, e che il dolce fosse il gusto preferito delle streghe, ma non ne era sicura''. 
\item[A.] Ah sì, giusto, non ci stavo ripensando. Perché è
   un'osservazione buttata lì, non legata alla storia. 
   Comunque, posso provare a dire come la vedo io per vedere se ho capito,
  visto che dei tre sono quella che ne sa di meno?
\item[N.] Certo, di' pure. 
\item[A.] Se ho visto che la stragrande maggioranza delle volte
  è stato esposto un cappello viola direi, anche senza il vostro Bayes,
  che dentro ci devono essere più streghe con il cappello di quel colore. 
  E quindi, per quello che abbiamo detto, non ci vuole una grande scienza,
  per capire che i cappelli devono essere quattordici viola e sette neri.
\item[G.] E quindi, come gusti?
\item[A.] Delle quattordici dal cappello viola ne avremmo dodici a cui piace il dolce e due a cui piace il salato, alle quali vanno aggiunte le sette streghe dal cappello nero. Quindi, in totale, a dodici piace il dolce e a nove il salato. Esattamente l'osservazione di Nora, ``\ldots\  e che il dolce fosse il gusto preferito delle streghe'', anche se non viene spiegato come ci arriva.
\item[N.] Forse semplicemente contando quante volte si erano arrabbiate quando era stato dato loro dolce, e quante volte nel caso opposto. 
\item[A.] Sì, può darsi, e per quello che la storia mi irrita, a parte le finezze matematiche. È possibile che nessuno si fosse accorto che quelle con il cappello nero preferivano solo salato, visto che avevano fatto tante prove, compreso preparare cibi né propriamente dolci né propriamente salati?
\item[N.] Diciamo che ``maggioranza di cappelli viola'' e ``gusto preferito dolce'' sono conclusioni consistenti, ma c'è un `ma': ``\ldots\  ma non ne era sicura''.
\item[A.] Allora non possiamo dire più niente? 
\item[N.] Beh, forse proprio niente no. Almeno ci potremmo limitare a dire quale delle due ipotesi ci sembra più credibile alla luce di quanto abbiamo osservato, ma non saprei da dove cominciare.
\item[G.] E è qui che entra in scena Bayes, o Laplace, se vogliamo essere più precisi\ldots\  ma lasciamo stare i dettagli storici. Insomma ci serve quello che oggigiorno, a torto o a ragione, va sotto il nome di ragionamento bayesiano. Mediante esso possiamo inferire la percentuale dei cappelli viola e da questa predire, seppur probabilisticamente, il colore del cappello del giorno dopo.
\item[N.] Scusa, ma adesso mi confondo. Ci stai dicendo che il ragionamento bayesiano entra in gioco per calcolare una cosa che nella fiaba non è calcolata, anche se nell'appendice c'è scritto che lo scopo della storia era proprio di illustrare il ragionamento bayesiano?
\item[G.] Proprio così. 
\item[N.] Dai, non è possibile, ti starai sbagliando. Le cinque autrici -- dico cinque! -- sono dell'ISTAT, il libro ha tutti i `timbri' dell'ISTAT e la postfazione è firmata da una professoressa, ordinario di statistica di Bologna. Poi arrivi tu, da Montorio, che nemmeno sei uno statistico, a dirci che nella favola non c'è nessun ragionamento bayesiano.
\item[G.] Eh, purtroppo è così, e la cosa sarebbe divertente, se non ci fosse da piangere.
\item[A.] Ma scusa, in alto all'ISTAT, qualcuno avrà pur\ldots 
\item[G.] In un paese come il nostro? Dove la ministra della pubblica istruzione è a dir poco scarsamente istruita --  mi riferisco al suo curriculum --  e commette strafalcioni di italiano, come \ldots
\item[A.] \ldots\ no, per carità, non ne parliamo\ldots
\item[N.] Va bene, ci fidiamo di te. Alba mi diceva di aver visto tuoi libri, anche in inglese, che avevano Bayes e bayesiano nel titolo. 
\item[G.] Troppo buone, come direbbe Fantozzi\ldots\ Ma a me non piace il principio di autorità, altrimenti vince l'ISTAT e la discussione è finita. E nemmeno quello di campanilismo, amicizia o parentela. Quindi facciamo così. Io vi spiego cosa si intende per ragionamento bayesiano e cosa dice il ``teorema che porta il suo nome'', mi riferisco a Bayes, come recita l'appendice. Poi giudicherete voi se c'è una cosa del genere nella storia. 
\item[A.] Spiegarci un teorema? Nel quale si sbagliano pure quelli dell'ISTAT? Chissà quanto deve essere complicato. No, senti, lasciamo stare.
\item[G.] Ah ah, complicato non lo è affatto, soprattutto se ci concentriamo sulla formula di aggiornamento del rapporto di probabilità. Tranquilla, non ci sono né quadrati né radici quadrate, e nemmeno logaritmi o funzioni trigonometriche.
\item[A.] Guarda in che pasticcio mi sono messa. Va be', prova, ma ti dico subito che se mi perdo lascio stare\ldots\ 
  Ma, scusate,  l'acqua per le tisane starà bollendo da un pezzo. 
  E mi sembra il momento buono pure per riprender fiato\ldots 
\end{description}
%
\begin{center}
{\Large  $[\,$ Pausa $\,]$ }
\end{center}  
\mbox{}\\
\begin{description}
\setlength\itemsep{0.5mm}
\item[A.] Eccomi. Thermos con acqua calda, 
the e tisane varie, fate voi. 
E ho portato pure qualche dolce 
per fare merenda, per rimanere in tema. 
Stamattina avevo preparato un ciambellone e una torta della nonna.
O forse Giulio preferisce qualcosa di salato, 
visto che si è presentato con un berretto nero\ldots 
\item[G.] Ah ah\ldots Grazie! Con the e tisane 
forse anche le streghe dal cappello nero
gradirebbero qualcosa di dolce\ldots \ 
Allora, affrontiamo finalmente in modo quantitativo questo benedetto ragiona\-mento bayesiano. Potete trovare varie formule su libri di testo e su internet, inclusi video su YouTube e, manco a dirlo, app per smart.
Ma la sostanza è la stessa. Io, come detto, vi mostro la mia preferita nel caso si abbiano due sole ipotesi, che nel nostro caso sono ``7 streghe con cappello viola'', che indichia\-mo
con $H_1$, e 
  ``14 streghe con cappello viola'', che indichiamo $H_2$.
  $H$ ricorda l'iniziale di `ipotesi' in latino, e anche in inglese.
\item[A.] Già, ``hypotheses non fingo'' \ldots
\item[G.] Qui non le dobbiamo fingere, perché le abbiamo dal contesto della fiaba.
\item[N.] Sì, ma non potremmo chiamare le due ipotesi $V_7$ e $V_{14}$, per non doverci ricordare chi è chi?
\item[G.] Ottima idea, tanto sono solo dei nomi, delle etichette, per così dire. 
\item[A.] E ora il teorema\ldots\ Mi sono già preparata mentalmente con apposita tisana\ldots 
\item[G.] Allora, cominciamo con gli `ingredienti' che ci servono -- per la formula, non per la tisana\ldots\ Innanzitutto ci serve il rapporto delle probabilità a priori delle due ipotesi.
\item[N.] Abbiamo detto che le riteniamo ugualmente possibili, le due ipotesi, quindi assegniamo loro probabilità 1/2 e 1/2, ossia 50\% e 50\%.
\item[G.] Ok, ma per quello che ci serve basta dire che sono uguali, in quanto siamo interes\-sati al loro rapporto, almeno per la formulazione che vi mostro fra un attimo.
\item[A.] Il rapporto di due cose uguali vale 1, o c'è qualche trucco nascosto?
\item[G.] Nessun trucco, allora scriviamo
    $$\frac{P(V_7)}{P(V_{14})} = 1\,.$$
\item[A.] E fin qui vi seguo ancora. Ma ci puoi dire subito qual è lo scopo di quello che stiamo facendo? O ripetercelo, perché forse ce l'hai già detto e mi sono persa. 
\item[G.] Vogliamo valutare come dobbiamo --  e insisto su `dobbiamo' perché l'aggiornamento va effettuato secondo una regola ben precisa e non alla buona --  \ldots\ dobbiamo cambiare il rapporto delle probabilità delle due ipotesi ogni volta che osserviamo un cappello all'ingresso della grotta. 
\item[A.] Insomma noi, abitanti di quel paesetto infestato dalle streghe, ogni giorno possia\-mo riaggiornare il rapporto delle probabilità, e quindi anche il valore di ciascuna probabilità, visto che la loro somma deve dare il 100\%. Dico bene?
\item[G.] Giusto. Ma per far ciò occorre valutare la probabilità di Nero e quella di Viola alla luce di ciascuna ipotesi, o, come si dice in gergo probabilistico, `condizio\-natamente a' ciascuna ipotesi. Se indichiamo `cappello nero' con $N$ e `cappello viola' con $V$, le probabilità che ci servono sono
  \begin{eqnarray*}
   && P(N\,|\,V_7)  \\
   && P(V\,|\,V_7)  \\
   && P(N\,|\,V_{14})  \\
  &&  P(V\,|\,V_{14})\,,
  \end{eqnarray*}  
  ove `$P()$'
  sta, banalmente, per `probabilità di' e la barra verticale si legge `condizionatamente a', o semplicemente `dato'. Il simbolo dentro le parentesi a sinistra della barra sta a indicare ciò di cui valutiamo la probabilità.
\item[A.] Quindi, se ho capito bene, $P(N\,|\,V_7)$ si legge ``probabilità di nero data l'ipotesi $V_7$'', ovvero ``probabilità di Nero condizionatamente all'ipotesi che ci siano sette streghe con il cappello viola'', e così via.
\item[N.] O più sinteticamente ``probabilità di $N$ dato $V_7$''.
\item[A.] Posso approfittare per fare pure io un diagramma come quelli
di prima. Mi stanno intrigando. Forse veramente li potrei usare a scuola anche per
questioni non prettamente quantitative. Ecco.
 \begin{center}
           \epsfig{file=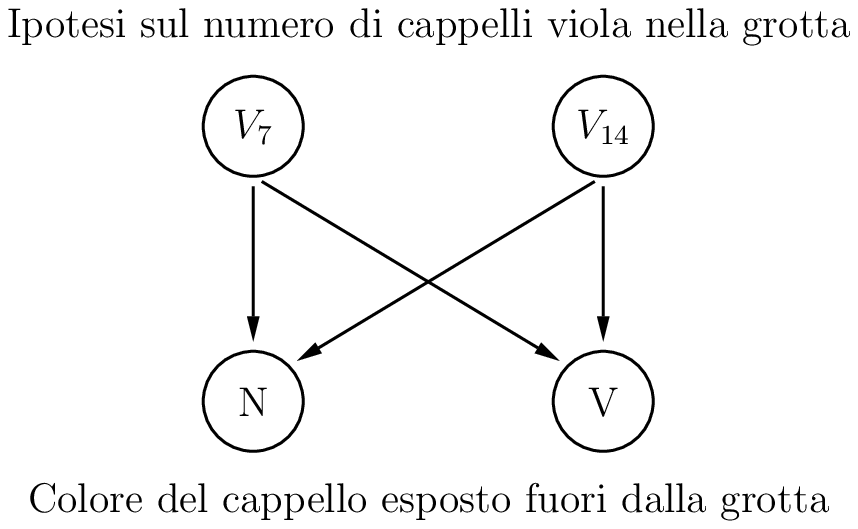,clip=,width=0.6\linewidth}
   \end{center}
\item[G.] Perfetto. Ci siamo capiti. Ora si tratta di valutare questi quattro numeri.
\item[N.] Facile. Comincerei da $P(V\,|\,V_7)$. Se ipotizziamo che le streghe con il cappello viola sono sette, allora
  $$P(V\,|\,V_7) = \frac{7}{21} = \frac{1}{3}\,$$
  e, per complemento, 
  $$P(N\,|\,V_7) = \frac{2}{3}\,.$$
  Se invece ne ipotizziamo quattordici abbiamo il contrario, ovvero
  \begin{eqnarray*}
    P(V\,|\,V_{14}) &=& \frac{2}{3}\\
    P(N\,|\,V_{14}) &=& \frac{1}{3}\,.
    \end{eqnarray*}
E ora possiamo mettere questi numeri sulla rete bayesiana che aveva disegnato
Alba e fare più spesse le linee che corrispondono alle probabilità
più elevate. 
\begin{center}
 \epsfig{file=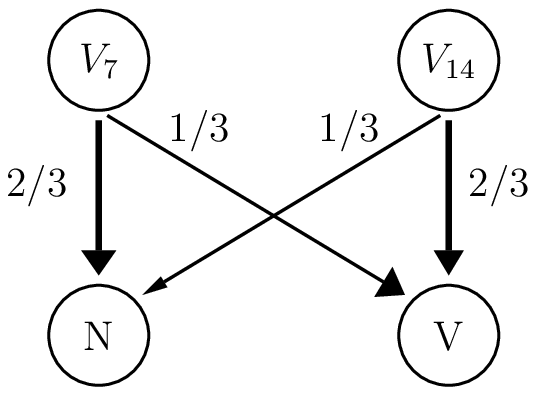,clip=,width=0.4\linewidth}
\end{center}
\mbox{}\vspace{-1.0cm}\mbox{}\\ 
  E ora?
\item[G.] Ora il tocco di Bayes! Una conseguenza del teorema che porta il suo nome ci insegna come aggiornare i rapporti di probabilità, secondo una regoletta facile facile. Nel nostro caso, dice semplicemente che se osserviamo un cappello di un certo colore {\em il fattore di aggiornamento è dato dal rapporto delle probabilità di osservare quel colore alla luce di ciascuna ipotesi}.
\item[N.] Ah, ecco cosa intendevi prima quando dicevi che si tratta di un vero fattore moltiplicativo. In effetti dipende solo da come le due ipotesi possano produrre più o meno facilmente tale colore e non da quanto esse sono probabili. 
\item[A.] Aspettate, mi sto perdendo, lo potete scrivere?
\item[N.] Posso provare io? Se ho capito bene, traducendo in formule quello che hai detto tu a parole, il fattore di aggiornamento vale, nel caso di Nero,
  $$   f  =  \frac{P(N\,|\,V_7)}{ P(N\,|\,V_{14}) } $$
  il cui valore è banalmente $$ \frac{2/3}{1/3} = 2\,.$$
Quindi la regola di aggiornamento dovrebbe essere
$$   \frac{ P(V_7)}{P(V_{14})} = \frac{  P(N\,|\,V_7)}{P(N\,|\,V_{14})} \times
\frac{P(V_7)}{ P(V_{14})}\,,$$
anche se mi vergogno di aver scritto una cosa del genere, e meno male che non c'è nessuno dei miei scolari a vedermi.
\item[A.] Che ci sarebbe di male?
\item[N.] Ci sarebbe, anzi c'è, il fatto che,
  affinché l'equazione sia valida,
 il primo rapporto a destra deve essere  unitario,
  mentre già sappiamo che vale due.
  C'è chiaramente qualcosa che non va, se la regola che Giulio ci ha detto è quella giusta.
\item[G.] Sì, la regola è giusta, ma bisogna fare un po' di attenzione alla notazione, perché $P(V_7)/P(V_{14})$ a destra e a sinistra dell'equazione non sono la stessa cosa.
\item[A.] Come non sono la stessa cosa? Pensavo di aver capito e ora mi stai confondendo nuovamente! Mi ci vuole un rinforzo di tisana.
\item[G.] Sì, scusa, Alba, ma non ti preoccupare, e se ho indotto Noemi a scrivere tale mostruosità mi devo essere espresso male io. Si tratta sempre di rapporti di quelle due probabilità, naturalmente, ma sotto condizioni diverse, e quindi i numeri, in genere, possono differire.
\item[A.] In che senso?
\item[G.] Facciamo un esempio pratico, con qualcosa che sicuramente hai in casa, anche se non ti chiederò di andare a prendere la scatola della tombola.
\item[A.] No, per carità! Ho appena messo via tutte le cose delle feste.
\item[G.] Non ti preoccupare, non serve averla fisicamente. Un po' di fantasia ce l'abbiamo. Quanto vale la probabilità di estrarre il `37' dal sacchetto dei numeri della tombola? Diciamo prima dell'estrazione.
\item[A.] Certo che prima. Che senso ha dopo?
\item[G.] Aspetta e vedrai.
\item[A.] Ci sono novanta numeri e quindi la probabilità di `37'
  vale $1/90$, giusto?
\item[G.] Ma ora immagina che io abbia già pescato il numero e ce l'abbia in mano, ma che non lo mostri a nessuno.
  Per la precisione non lo vedo nemmeno io.
  Quanto vale la probabilità che quel numero sia `37'?
\item[A.] Ancora 1/90. Perché mai dovrebbe cambiare? Anche se l'estrazione è già avvenuta. Quello che è importante è conoscerne o meno il risultato. Ah, ora capisco cosa intendevi dire
  quando precisavi ``prima dell'estrazione''!
\item[N.] E difatti non cambia, giusto? Per quale motivo razionale dovresti cambiare idea? A meno che non ti accorgi che anche prima avevi sbagliato --  ma questa è un'altra storia.
\item[G.] Ma ora immagina che io permetta a Noemi di sbirciare
  il numero, e poi
  le chieda di dirti se è pari o dispari, niente di più.
  E che te lo dica sottovoce all'orecchio, in modo che io non possa
sentire.   Cambia qualcosa?
\item[A.] E adesso che c'entra che me lo deve dire sottovoce?
\item[G.] Aspetta e vedrai.  
\item[A.] Se mi dice che è pari, la situazione cambia drasticamente, perché sicuramente non è il `37'. Ma anche se mi dicesse dispari, la probabilità cambierebbe, perché ora c'è solo una possibilità su quarantacinque, ossia una probabilità di $1/45$. 
\item[G.] Vedi? Si tratta in effetti di tre probabilità di `37' alla luce di tre stati di informazione diversi. Il primo è quello iniziale, che scriviamo $P_0(`37')$, e che vale 1/90. A pedice di P ci abbiamo messo `0' per ricordare che si tratta di quella iniziale, insomma ``all'istante zero'', come si dice in Fisica. Se il numero è stato estratto ma tu non hai alcuna informazione specifica rimane sempre $P_0(\mbox{`37'}) = 1/90$.
\item[A.] Qualsiasi numero sia, lui è e lui rimane.
  Insomma, se capisco bene, stai insistendo sul
  fatto che il nostro giudizio probabilistico non
  è legato  al fatto che possa uscire un altro numero al posto suo,
  ma al nostro stato di conoscenza.
\item[G.] Proprio così! Se ti è stato detto che il numero è pari
  la probabilità diventa
  $$P(\mbox{`37'}\,|\,\mbox{`pari'}) = 0\,,$$
  mentre se ti è stato detto che è dispari vale
  $$P(\mbox{`37'}\,|\,\mbox{`dispari'}) = \frac{1}{45}\,.$$
\item[A.] Mi sembra di aver capito. Il condizionamento cambia il valore di probabilità e quindi è meglio usare simboli diversi per capire di cosa si sta parlando.
\item[G.] Bene, vedo che ci siamo capiti. 
\item[N.] Scusa, Giulio, possiamo usare la tua regola\ldots
\item[G.] Sì, mia, magari!
\item[N.] Intendo la regola di Bayes con i rapporti, per vedere se le cose tornano?
\item[G.] Ottima idea! Vuoi provare?
\item[N.] Certo, sembra facile, se consideriamo solo due possibilità: `37' o `non-37', ossia uno qualsiasi degli altri numeri. E stavolta utilizzerei la notazione $H_1$ e $H_2$ per indicare le due ipotesi.
  Quindi ci servono le probabilità di `pari' e di `dispari' dato $H_1$,
  che valgono banalmente  0 e 1, rispettivamente:
  \begin{eqnarray*}
    P(\mbox{`pari'}\,|\,H_1) &=& 0 \\
    P(\mbox{`dispari'}\,|\,H_1) &=& 1\,.
    \end{eqnarray*}
\item[A.] E fin qui ci sono\ldots
\item[N.] Poi ancora $P(\mbox{`pari'}\,|\,H_2)$, che vale 45 su 89, e $P(\mbox{`dispari'}\,|\,H_2)$, che vale 44 su 89.
\item[A.] Aspetta, perché 45 nel caso `pari' e 44 nel caso `dispari'?
\item[N.] Perché $H_2$ comprende i numeri da 1 a 36 e da 38 a 90,
  e quindi l'insieme è costituito da 45 pari e 44 dispari. Riscriviamo in bell'ordine
    \begin{eqnarray*}
    P(\mbox{`pari'}\,|\,H_2) &=& {45}/{89} \\
    P(\mbox{`dispari'}\,|\,H_2) &=&   {44}/{89} \,.
    \end{eqnarray*}
\item[A.] Ah grazie, era ovvio, scusa. Quindi se tu hai visto '37'
  mi devi dire per forza dispari; se invece hai visto uno degli
  altri numeri è più facile che mi dici `pari' che `dispari', anche se di poco:
  45/89 rispetto a 44/89. Fin qui ci sono.
\item[G.] Poi servono le probabilità iniziali, anzi concentriamoci direttamente sul rapporto delle probabilità iniziali. 
\item[N.] Facile: avendo una possibilità contro le rimanenti 89, il rapporto è
            $$\frac{P_0(H_1)}{ P_0(H_2)} = \frac{1}{89}\,.$$
Mentre invece quella specie di fattore di aggiornamento, ossia il rapporto delle probabilità che io possa dire a Alba `pari' alla luce delle due ipotesi, vale
            $$\frac{P(\mbox{`pari'}\,|\,H_1)}{ P(\mbox{`pari'}\,|\,H_2)} = \frac{0}{45/89}\,,$$
che è banalmente zero. E per qualsiasi numero lo moltiplichiamo il risultato sarà sempre zero. Quindi se le dico `pari' la probabilità che il numero possa essere $H_1$,
ovvero `37', sarà nulla.
\item[A.] Bella scoperta! Ci volevano tutti questi conti per dirlo?
\item[G.] In effetti questo caso è talmente facile che, come si dice, stiamo rompendo una noce con un maglio. Stiamo solo controllando che la regoletta che discende dal teorema di Bayes dia il risultato che sappiamo essere giusto.
\item[A.] Ok, non volevo essere polemica. 
\item[N.] E ora passiamo al secondo caso, ovvero quando dico a Alba che è dispari. In questo caso abbiamo
             $$\frac{P(\mbox{`dispari'}\,|\,H_1) }{ P(\mbox{`dispari'}\,|\,H_2)} = \frac{1}{44/89} =  \frac{89}{44}\,,$$
che è quasi due, e quindi il rapporto delle due probabilità praticamente raddoppia.
\item[A.] Quant'era all'inizio? 
\item[N.] Era 1 a 89, in quanto una possibilità contro le ottantanove rimanenti. L'avevamo già scritto. 
\item[A.] Ah, sì, scusa, giusto!
\item[N.] E questo numero va moltiplicato per il fattore di aggiornamento, dato dall'altro rapporto, che abbiamo visto essere $89/44$.
\item[G.] A proposito, sapete il nome che viene dato a questo fattore di aggiornamento dei rapporti di probabilità?
\item[N.] Nessuna idea.
\item[G.] Fattore di Bayes.
\item[A.] Ci mancava\ldots
\item[G.] Anche se qualcuno dice che andrebbe chiamato fattore di Bayes-Turing. 
\item[A.] Turing, quello del codice Enigma? Ho visto il film, {\em The imitation Game}, ma non parlavano di questo fattore.
\item[G.] Sì, lui, un grande! Brutta fine, lasciamo stare\ldots\ E dal film non si capisce un granché, a parte che era un tipo geniale, e i suoi drammi personali.
\item[N.] E perché questo fattore dovrebbe portare anche il suo nome?
\item[G.] Perché secondo un suo collaboratore a Bletchley Park, un certo I.J. Good, il metodo di Turing per capire  --  direi inferire  --  la chiave giornaliera del codice Enigma era basata proprio su fattori di aggiornamento del genere.
\item[A.] Toh, vallo a pensare! Immagino un problema ben più difficile di quello di calcolare la probabilità che il numero sia `37' sapendo che è pari\ldots\ e senza i nostri computer. Mi ricordo nel film quella specie di grossi armadi con fili e ingranaggi.
\item[N.] Scusa, come fa di nome questo Good?
\item[G.] Sai che non lo so? Mentre Turing è Alan Turing, di Good ho sempre sentito parlare come `I.J. Good', come fosse un tutt'uno, forse perché `Good' da solo si confonderebbe. Mi son fatto un'idea che lo chiamasse così anche la madre da bambino\ldots \ 
Insomma, I.J. Good chiama così questo fattore in un breve articolo segretato ai suoi tempi e reso pubblico solo pochi anni fa~\cite{IJGood}.
\item[A.] Interessante! Poi ci dai delle informazioni in proposito, soprattutto sul ruolo di Turing nell'impresa per decriptare i messaggi dei tedeschi.
\item[G.] Sì, guarda, senza andare nel dettaglio, vi raccomando la conferenza che
  una giorna\-lista di cui non ricordo mai il nome
  ha tenuto anni fa alla sede di Google per presentare il suo libro
  {\em The theory that would not die}~\cite{McGrayne} --
  il titolo del libro invece lo ricordo. Trovate facilmente il video della conferenza su
  Youtube.\footnote{\,\url{https://www.youtube.com/watch?v=8oD6eBkjF9o&t=11s}}\\
  Ma torniamo al nostro problema, anche perché capite bene che, se quella signora ci ha scritto un tomo di varie centinaia di pagine, sull'argomento ci sarebbe un bel po' da dire.
\item[N.] Allora vado avanti. Il rapporto iniziale vale $1/89$, abbiamo detto. Il fattore di aggiornamento, questo fattore di Bayes\ldots
\item[A.] Bayes-Turing mi piace di più, adesso che sto collegando questa formuletta a eventi talmente importanti.
\item[G.] Anche a me, se proprio gli dobbiamo dare un nome, con l'acronimo BTF,
  {\em Bayes-Turing Factor}, in inglese.
\item[N.] Riepilogando,  nel nostro caso
  il fattore di Bayes-Turing è 89/44 e il rapporto finale vale quindi
$$ \frac{89}{44}  \times  \frac{1}{89} = \frac{1}{44}\,.$$ 
Il risultato è quindi
    $$ \frac{P(H_1\,|\,\mbox{`dispari'}) }{ P(H_2\,|\,\mbox{`dispari'})} = \frac{1}{44}\,. $$ 
\item[A.] Quindi, fallo dire a me per vedere se ho capito:
  l'ipotesi $H_2$, ovvero `non-37', era inizialmente 89 volte 
più probabile dell'altra;
  ora, alla luce del fatto che io so che il numero estratto è dispari,
  è diventata soltanto 44 volte più probabile.  
\item[N.] E quindi possiamo dire che ora la probabilità di '37' è
  pari a $1/45$, esattamente il numero a cui eravamo arrivati prima.
\item[G.] Ma in realtà non  arriviamo tutti e tre a questo
  valore di probabilità.
\item[N.] E perché no?
\item[G.] Perché nell'ipotesi di lavoro, per quanto fittizia,
  non abbiamo lo stesso stato di informazione:
  \begin{itemize}   
  \item Alba ha solo sentito `dispari' e quindi 
    lei \underline{deve} riaggiornare
    la probabilità a $1/45$;
  \item tu hai visto il numero e quindi \underline{sai} qual è con certezza;
  \item io infine non ho visto né sentito niente e quindi \underline{resto}
    della mia opinione. 
  \end{itemize}
\item[A.] Ecco il mistero di farmi sussurrare all'orecchio
  se era pari o dispari!
  Comunque abbiamo riottenuto, almeno per me, il risultato precedente.
\item[G.] Giusto. Ma ora ora capite bene che possediamo un potente strumento matema\-tico per aggiornare i rapporti di probabilità, e possiamo usarlo anche quanto l'intuizione non ci aiuta, o addirittura ci condurrebbe a conclusioni sbagliate. 
\item[A.] E quindi il famoso ragionamento bayesiano è tutto lì!?
\item[G.] Sì, almeno nel caso di due ipotesi. 
E anche se sono tante esso vale per ciascuna coppia di ipotesi. 
\item[N.] Quindi, se sappiamo come il nuovo rapporto di probabilità di ogni coppia di ipotesi possiamo normalizzare e trovare la probabilità di ogni ipotesi. 
\item[A.] Sì, va be', normalizzate pure, io mi accontento di due possibilità.
\item[G.] Dai, non ti spaventare, sono soltanto dettagli tecnici, se li vogliamo chiamare così. La sostanza è quella che abbiamo visto: {\em il rapporto di probabilità viene aggiornato moltiplicandolo per il fattore di Bayes-Turing, che è dato dal rapporto delle probabilità con cui l'informazione potrebbe scaturire da ciascuna ipotesi.} Bene, ora possiamo andare avanti a risolvere il problema inferenziale di capire quante streghe di ciascun gruppo ci sono nella caverna in base ai colori dei cappelli osservati ogni giorno.
\item[N.] Anche in questo caso il calcolo dei due fattori è facile. Il primo
  l'avevamo già calcolato. Eccolo lì su quel foglio:
  $$\frac{P(N\,|\,V_7)}{P(N\,|\,V_{14})} = \frac{2/3}{1/3} = 2\,,$$\\
  e ci scriviamo affianco anche il secondo, che è 
  $$\frac{P(V\,|\,V_7)}{P(V\,|\,V_{14})} = \frac{1/3}{2/3} = \frac{1}{2}\,,$$
\item[A.] Quindi, fatemi vedere se ho capito: se si osserva Nero la probabilità di $V_7$ si raddoppia rispetto a quella di $V_{14}$, mentre se si osserva Viola succede l'opposto, raddoppia la probabilità di $V_{14}$ rispetto a
  quella di $V_7$.
\item[N.] Giusto, perché è più facile che $V_7$ dia Nero, mentre è più facile che $V_{14}$ dia Viola. E il giorno seguente?
\item[G.] Si va avanti nello stesso modo. 
\item[A.] Quindi, se per ipotesi davanti alla caverna comparisse Nero per quattro giorni di fila? Perché no\ldots
\item[G.] E anche 10, 20 o 100 giorni, seppure tale evento è a priori poco probabile, ma non impossibile. E una volta che si è verificato ce lo teniamo. 
\item[N.] Facciamo il caso di quattro cappelli neri 
quattro giorni consecutivi. 
\item[A.] E, se permettete, provo a fare il diagramma del modello. Eccolo:
\begin{center}
 \epsfig{file=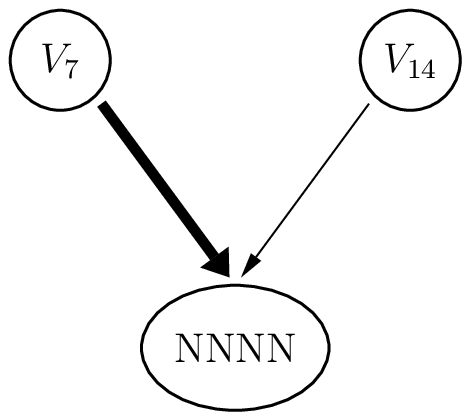,clip=,width=0.34\linewidth}
   \end{center}
Ho anche fatto la freccia a sinistra bella cicciotta, 
perché mi pare ovvio che $V_7$ è quella che dà più facilmente 
quattro neri su quattro.  Sembra quasi troppo semplice, spero sia giusto.
\item[G.] Tranquilla, hai un futuro di designer di reti bayesiane\ldots\ 
Sì, a rigore, bisognerebbe mettere tutte le altre sequenze
di quattro estrazioni, ma sono irrilevanti, in quanto 
quello che conta è il condizionamento dovuto a $NNNN$.
\item[A.] Giusto, ci sarebbero anche le altre sequenze, ecco
  \begin{center}
 \epsfig{file=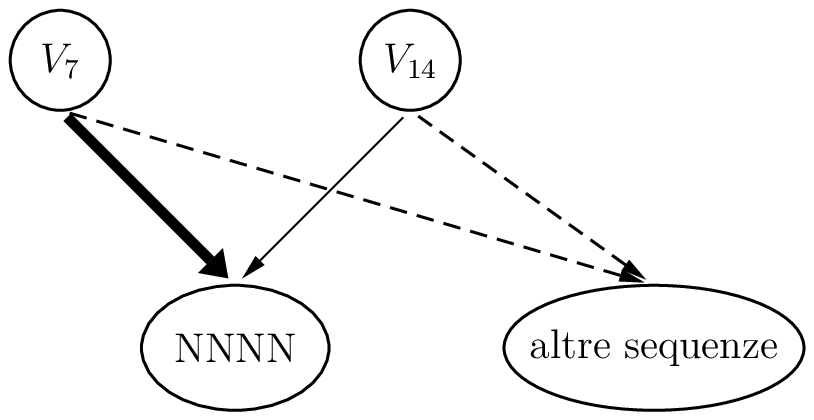,clip=,width=0.59\linewidth}
  \end{center}
  in cui ci tracciato linee tratteggiate per non stare troppo a pensare
  a quali possono essere le probabilità.
\item[N.] Beh, non sono complicate:
  \begin{itemize}
  \item quelle che partono da $V_7$ valgono
  $(2/3)^4 = 16/81$ e, per complemento, $65/81$;
  \item quelle che partono da $V_{14}$ sono invece 
    $(1/3)^4 = 1/81$ e $80/81$.
  \end{itemize}  
\item[A.] Beata te!
\item[G.] E invece le probabilità
  delle ipotesi alla luce della sequenza osservata?
\item[N.] Facile. Il rapporto delle probabilità si raddoppia 4 volte: 2, 4, 8, 16, o in modo compatto $2^4$. L'ipotesi $V_7$ diventa quindi 16 volte più probabile di $V_{14}$, giusto? 
\item[G.] Giusto, e scriveremo allora
         $$  \frac{ P(V_7\,|\,4N) }{ P(V_{14}\,|\,4N)} = \frac{16}{1},$$
  ove, con `$|\,4N$' abbiamo reso esplicito che tale valutazione è condizionata
  dalla osservazione di quattro cappelli neri. Otteniamo quindi
          $$ P(V_7\,|\,4N)  = \frac{16}{17} $$
e
          $$ P(V_{14}\,|\,4N) = \frac{1}{17}\,. $$
\item[A.] Scusa com'è che sei passato dal rapporto 16 a 1 alle probabilità 16/17 e 1/17.
\item[N.] Questo è facile, basta rinormalizzare numeratore e denominatore alla loro somma, come abbiamo fatto prima con 1/44, quando eravamo interessati a `37' estratto dal sacchetto della tombola sapendo che si tratta di un numero dispari. Un dettaglio aritmetico. 
\item[A.] Sì, grazie! È che prima nemmeno mi ero posta il problema, perché sapevo già come andava a finire. E mi sembra pure di aver capito cosa intendete per `normalizzazione': scrivere le probabilità individuali in modo tale che la somma dia 1. 
Infine, fatemi vedere se ho capito, se successivamente osserviamo 4 volte il cappello viola il fattore di Bayes-Turing si riduce di un fattore 16, e torniamo al punto di partenza.
\item[G.] Sì, perché il fattore di aggiornamento globale è pari al prodotto di quelli parziali. Ne segue che quattro cappelli neri e quattro viola producono un fattore di aggiornamento globale di
  $16 \times (1/16)$, che è uguale a uno. E notate pure come l'ordine delle osservazioni non conta.
\item[N.] In pratica quindi quello che conta è solo la differenza delle volte che abbiamo visto i due colori. Non solo indipendentemente dall'ordine, ma anche dalla particolare sequenza, ad esempio $NNVVVVVNVNV$ ci dà lo stesso fattore di aggiornamento di $VNVNVVVVVNN$. Interessante!
\item[G.] In questo caso sì, perché `cappello nero' raddoppia
  il rapporto definito in quel modo, mentre `cappello viola' lo  dimezza.
  Quindi il fattore di aggiornamento vale $2^{\#N}/2^{\#V}$, ovvero $2^{\#N - \#V}$, avendo indicato con $\#N$ il numero di volte che abbiamo visto il cappello nero, e similmente per il viola.
\item[N.] Quindi se in altro problema l'effetto $A$ avesse triplicato il fattore di Bayes-Turing e l'effetto $B$ lo avesse dimezzato, il risultato non sarebbe dipeso semplicemente dalla differenza, essendo il fattore di aggiornamento pari a $3^{\#A}/2^{\#B}$.
\item[G.] Ottima osservazione, grazie.
\item[N.] Tutto qui? In effetti la matematica è veramente elementare. Il teorema di Pitagora, con quadrati e radici quadrate, è ben più complicato. Sì, qui possiamo avere delle potenze, ma solo se vogliamo scrivere in modo compatto il fattore di aggiornamento globale, equivalente a più aggiornamenti in successione.
\item[G.] Assolutamente!
\item[N.] E, se permettete, a parte le questioni numeriche, importanti se vogliamo giungere a risultati quantitativi, mi sembra che il succo di questa inferenza bayesiana, nel problema specifico, consista nel fare affermazioni probabilistiche su quello che non possiamo vedere, ovvero il numero di streghe appartenenti a ciascun gruppo, alla luce di quello che ci è concesso di vedere, ovvero la sequenza dei colori dei cappelli posti davanti alla grotta. E per far ciò è fondamentale `sapere', o insomma ipotizzare in qualche modo, come ciò che è osservabile empiricamente è logicamente legato, seppur probabilisticamente, alle ipotesi che si vogliono inferire.
  Insomma, dobbiamo aver definito un modello causale probabilistico, 
come abbiamo imparato.
  E se poi lo rappresentiamo graficamente ancora meglio, perché,
  come spesso succede, le rappresentazioni grafiche spesso aiutano
la comprensione e quindi la soluzione di problemi. 
\item[A.] E, scusatemi se filosofeggio, forse perché anche qui abbiamo una caverna, ma questo tipo di ragionamento mi fa venire in mente il famoso mito di Platone. Non ci è dato di guardare dentro ma possiamo vedere solo quello che si manifesta fuori, che però deve essere connesso in qualche modo a quello che succede dentro.
\item[G.] Sì, buona idea, me lo dovrei rivedere il mito della caverna.
\item[A.] Sarebbe divertente scoprire che Platone avesse già effettuato ragionamenti bayesia\-ni
  {\em ante litteram}.
\item[G.] Sì, divertente, ma non ci credo molto, chissà\ldots\  Ma, tornando a noi, la questione non è se c'è un ragionamento bayesiano, seppur primordiale, negli scritti di Platone ma, più terra terra, se c'è nella fiaba delle streghe di Bayes.
\item[N.] Non che mi sembri.
\item[G.] Tranquilla, la risposta è no. Nella fiaba non c'è traccia di ragionamento bayesiano, punto e basta. Direi che forse è più facile trovarlo in Platone, leggendo fra le righe\ldots\  chissà\ldots
\item[A.] Ma -- e scusa se faccio l'avvocato del diavolo  --  non è che il ragionamento bayesiano c'entra quando estraggono a sorte uno dei medaglioni con S e D per decidere cosa dare alle streghe?
\item[G.] Ah ecco, hai toccato un altro punto dolente della storia. Tu Noemi che ne pensi? Oops, scusate, deve essere il mio\ldots\  sì,
  un attimo che rispondo\ldots 
\item[N.] Approfitto pure io per avvertire che la cosa si sta facendo un po' lunga.
  Fortuna che siamo invitati da amici e non devo preparare la cena.
  Io avevo pensato che avremmo fatto una analisi del testo, che è di cinque o sei
  pagine. Invece, ridendo e scherzando abbiamo spaziato da Platone a Fantozzi,
  passando per Leibnitz e David Hume\ldots
\end{description}
%
\begin{center}
{\Large  $[\,$ Pausa $\,]$ }
\end{center}  
\mbox{}\\
\begin{description}
\setlength\itemsep{0.5mm}
\item[G.] Ecco, scusate, era mio fratello. Voleva sapere per stasera. 
  Avevo silenziato lo smart ma non il telefonino\ldots\  Non mi guardate così, ho un telefonino per telefonare e uno smart per\ldots\  smartare. Il telefonino è un vecchio Nokia,
  che se casca rimbalza senza che gli succeda niente e
  le batterie durano quasi una settimana\ldots\
  e mentre telefoni puoi usare lo smart per cercare informazioni, e così via.
Ma torniamo a noi. Quindi, ti chiedevo, che ne pensi della strategia decisionale di Nora?
\item[N.] Io quando l'ho letta ho pensato che se questo è quello che scrivono quelli del'ISTAT devo essere proprio stupida, perché quella strategia non
la capisco priprio. E tu mi stai dicendo che forse ho ragione io?
\item[G.] Non so ancora quello che pensi, ma siccome una persona di buon senso non farebbe mai quello che suggerisce Nora, che qui rappresenta l'ISTAT, sono pronto a scommettere che hai ragione. Ci spieghi cosa avresti fatto tu, ossia perché la strategia dell'estrazione del medaglione non ti convince?
\item[N.] Facile. Io ho ragionato così, e dimmi se sbaglio. Innanzitutto prendiamo per buona l'interpretazione della filastrocca: cappelli neri vogliono solo salato; delle streghe con cappello viola solo una su sette preferisce salato e le restanti dolce. Giusto?
\item[G.] Sì, su questo eravamo già d'accordo e non possiamo farci niente, essendo un assunto della storia.
\item[N.] Se quindi vedo il cappello nero offro sempre cibo salato, e non sbaglio mai.
\item[A.] E per questo non bisogna fare grandi studi di probabilità e statistica, solo quegli stupidi di bayesiani potevano star tanto a lambiccarsi il cervello. Bastava fare una prova e continuare con quello che aveva funzionato.
\item[N.] Se invece vedo il cappello viola, porterei sempre dolce.
  In media sei volte su sette mi va bene e solo una volta su sette dovrò subire gli scherzi delle streghe\ldots
\item[A.] A proposito, scusate, ma questa è un'altra cosa che non fila nella storia, a parte i dettagli matematici. Si dice che le ventuno streghe si installano nella caverna e che, siccome sono pigre, non vogliono cucinare. 
\item[N.] E perché no? Questa è una fiaba e la dobbiamo accettare.
\item[A.] Sì. Va bene, ma c'è un limite a tutto! È che sembra che mangi solo la strega che è stata estratta, e per di più solo se è fortunata e le danno quello che gradisce. Insomma la maggior parte di loro digiuna, e qualcuna anche a lungo, se non viene estratta per molto tempo.
  E magari la volta che è estratta i bayesiani le danno pure il cibo che a lei non piace. Ecco perché si in\ldots\  si inquieta così tanto,
  per usare il verbo che avrebbe accettato
  il mitico maestro Guido, che ci segnava in rosso anche `arrabbiarsi'. 
\item[G.] Ah ah, a questo non ci avevo pensato, concentrato com'ero sulle questioni pro\-babilistiche. Ma lasciamo che Noemi ci dica come farebbe lei.
\item[N.] Quello che farei io l'ho detto e lo ripeto: cappello nero, sempre salato; cappello viola, sempre dolce. E così minimizzo le conseguenze negative. E loro invece che ti fanno? Intendo quelli dell'ISTAT, nei panni della figlia del fornaio. Nel caso di cappello viola estraggono a caso un medaglione su sette, sei dei quali hanno inciso D e uno S, per `dolce' e `salato'. Ma è chiaro che così si peggiorano le cose: il medaglione e la strega estratti sono eventi indipendenti e bisogna moltiplicare le probabilità\ldots
\item[A.] Scusa, Noemi, ma fatico a seguirti. Riesci a spiegarlo in modo semplice, eventualmente con qualche disegno?
\item[N.] Ok. Ma tanto per capirci, sei d'accordo che se fai come dico io, ti va sempre bene quando vedi il cappello nero, mentre in caso di cappello viola ti va male, in media,
  una volta su 7? E una volta su 7 è circa il 14\%.
\item[A.] Questo è chiaro. Non ho invece capito se con il metodo di Nora o, come dici tu, dell'ISTAT, ci si guadagna o ci si perde rispetto a quello che proponi tu. Quando avevo letto la storia, fidandomi delle autrici, che sono ben cinque, e considerando il fatto che la strategia era cervellotica, ho pensato che un vantaggio lo dovesse almeno avere. Cioè avevo assunto che fosse la migliore strategia possibile e che quindi ci si dovesse guadagnare in termini di quieto vivere nel villaggio. Tu invece mi stai dicendo, e pure Giulio sembra d'accordo, che la strategia è meno performante, come si dice oggi, di quella che sarebbe venuta in mente a qualsiasi bambino delle elementari. Insomma, quelli dell'ISTAT sembra che non ne azzecchino una. Spero almeno che sappiano cosa dicono per questioni ben più serie di questa fiaba, come quando parlano di povertà e disoccupazione nel paese, di vita attesa e pensioni, e di tante altre cose.
\item[N.] Su questo non so giudicare. So solo che se uno sbaglia a fare calcoli elementari e poi usa quei ragionamenti per decidere se lanciare bombe atomiche sarei molto preoccupata. Ma che dico lanciare bombe! Basta pensare ai problemi delle vaccinazioni. Va be', incrociamo le dita e torniamo alle nostre streghe. Allora, ragiona, la strega è estratta a caso, e di quelle dal cappello viola solo una su sette gradisce il salato. Quindi abbiamo sette possibilità, che possiamo associare alle colonne di una scacchiera. Ah eccone là una! La posso prendere?
\item[A.] Ma non sappiamo se le streghe dal cappello viola sono sette o quattordici.
\item[N.] Non fa niente. Ragiona per il momento come se fossero sette. Tanto quello che conta sono le proporzioni. Ecco, questa è la scacchiera, che ha pure impressi numeri e lettere per identificare le caselle. Innanzitutto ignoriamo l'ottava colonna, la H, mentre ciascuna colonna da  A a F la associamo a una strega dal cappello viola a cui piace il dolce. Infine la settima colonna, G, la associamo alla strega dal cappello viola che preferisce salato.
\item[A.] Ok, e ho pure capito che per questo ragionamento sette o quattordici, ovvero $V_7$ o $V_{14}$, è la stessa cosa. È come se raddoppiassimo le colonne, sovrapponendone altre a quelle proprie della scacchiera.
\item[N.] Le righe da 1 a 6 le associamo invece ai medaglioni con la scritta D, dolce; la riga 7 al medaglione con su scritto S, salato, e la riga 8 la ignoriamo.
\item[A.] Forse ho capito, posso seguitare a dirlo io?
  Correggetemi se sbaglio.
  L'estrazione della strega e quella del medaglione corrispondono a una casella della scacchiera nella regione ammessa, ovvero abbiamo 49 possibilità equiprobabili, giusto? E quand'è che la strega è contenta? Nella nostra analogia è quando le colonne sono da A a F, mentre le righe sono da 1 a 6; oppure quando si estrae la casella G7. Abbiamo quindi 37 casi su 49. Nei rimanenti 12 casi su 49 la strega si arrabbia. Giusto?
\\ \mbox{} \vspace{-0.5cm} \\ \vspace{-0.5cm}  
\begin{center}
\epsfig{file=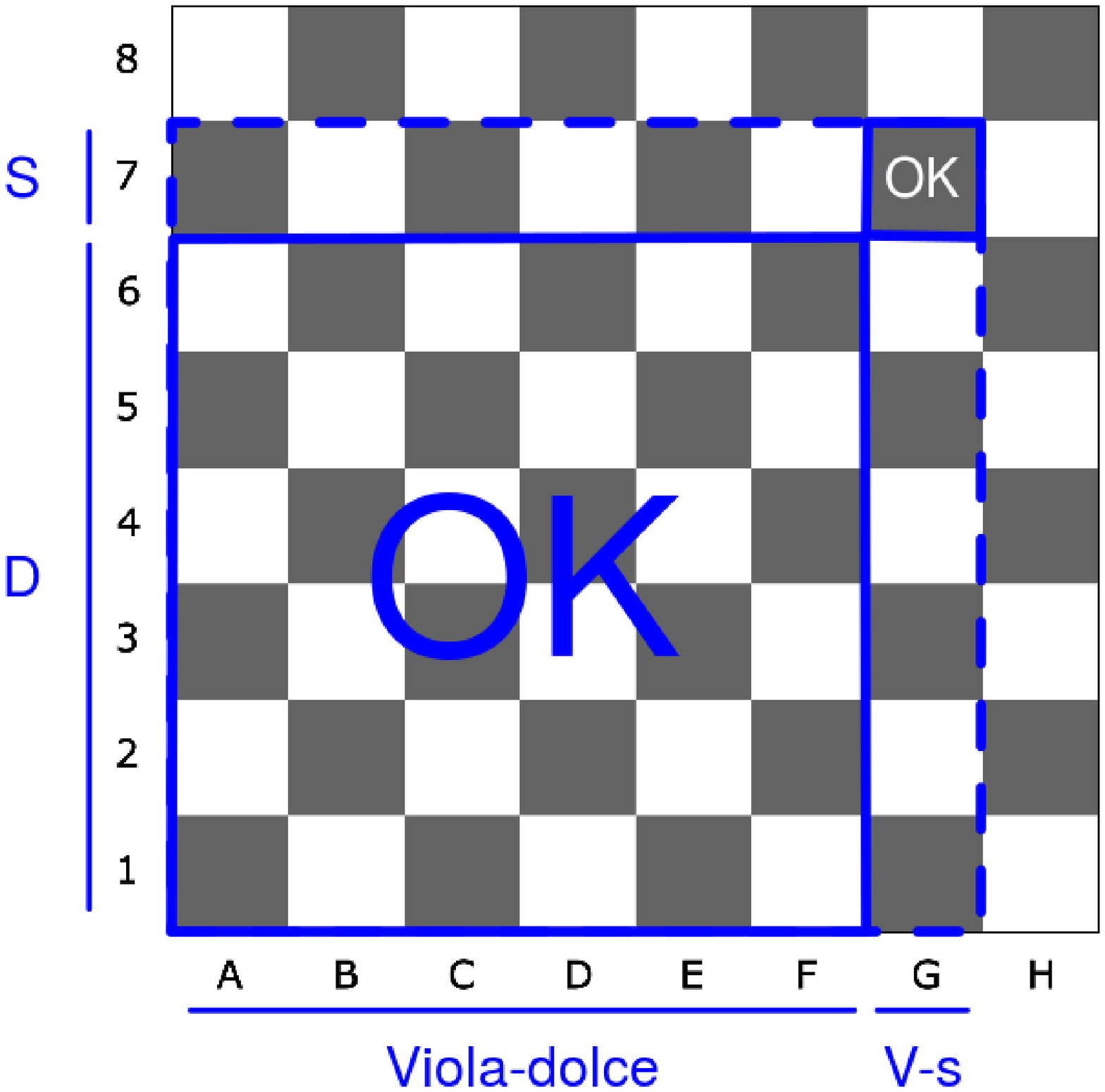,clip=,width=0.55\linewidth}
\end{center}
\item[N.] Esattamente quello che avevo calcolato io. Mentre con la strategia `a intuito' la probabilità che la strega si arrabbi è 1/7, ovvero 7/49, se moltiplichiamo numeratore e denominatore per 7.
\item[A.] Perché ti devi complicare la vita? 1/7 non andava bene?
\item[N.] In realtà ce la semplifichiamo perché dal ragionamento che hai appena fatto è chiaro che tutte le altre probabilità verranno espresse in quarantanovesimi. Insomma, con la strategia che adotterebbe un bambino delle elementari la proba\-bilità che la strega si arrabbi vale 7/49, corrispondente alle sette caselle
  consentite dalla colonna G rispetto alle quarantanove possibili. 
  Con la strategia ISTAT essa aumenta invece a 12/49. Quasi raddoppia!
 Hai capito come mai non riuscivo a capire dove sbagliavo?
\item[G.] Tranquilla, non sei tu a aver sbagliato, purtroppo\ldots
\item[N.] Purtroppo come?
\item[A.] Non ti offendere, Noemi. Giulio intende dire che 
sarebbe stato meglio che ti fossi sbagliata tu piuttosto che tutti quei cervelloni dell'ISTAT, fra autrici e tutti coloro che il libro lo hanno approvato, letto e pubblicizzato. E poi ci dobbiamo fidare di quello che ci raccontano sulle pensioni e su altro? Spero almeno che quanti polli mangiamo a testa, in media, lo sappiano calcolare ancora! Senza considerare che, quando le cose diventano più importanti, oltre alla possibilità di sbagliare per incompetenza o distrazione, mettici pure le pressioni di politici e dei gruppi di potere\ldots\  Va be', scusate, la smetto.
\item[G.] Io pure sono un po' preoccupato, ma tendo a essere fiducioso. Ci saranno sicuramente tante brave persone, oneste e competenti. Insomma, è come se qualcuno esterno al mondo della scuola dovesse giudicare il corpo docente dagli strafalcioni della ministra.
\item[N.] E magari il presidente dell'ISTAT, se proprio vogliamo parlare dei vertici, il libro non lo ha nemmeno letto.
\item[G.] Prima della pubblicazione forse no, non saprei cosa dire.
  Ma dopo potrebbe averlo letto, anzi avrebbe dovuto leggerlo,
  o almeno farlo leggere, se è un dirigente serio, a qualche persona di fiducia,
  se lui non aveva tempo.
\item[A.] Cosa vorresti dire?
\item[G.] Dico questo perché un mio collega, ora in pensione, che di queste cose ne capisce e a cui avevo segnalato il libro, gli ha addirittura scritto, tanto era rimasto scandalizzato da quello che vi aveva letto.
\item[A.] E che gli ha risposto?
\item[G.] Magari gli avesse risposto! Capisci come siamo conciati? E non è il primo che capita. Lo scorso anno ha addirittura pubblicato un libro sull'argomento in cui, fra l'altro, si mettevano in evidenza errori che si compiono, a diversi livelli, dovuti a ignoranza su probabilità e
  statistica\,\cite{Pallottino}.
  E basta cercare il suo nome su Google e ottieni in una frazione di secondo informazioni su di lui e sul libro che ha scritto.
  Insomma a uno come lui non puoi non rispondere! Al più gli fai
  rispondere da qualcuno di fiducia. 
  E comunque, chiunque sia, come minimo si valuta se quello
  che scrive è degno di essere preso in considerazione.
\item[A.] Uhm\ldots\  non mi sembra un comportamento serio, a meno che di queste cose non ci capiscano niente manco loro, intendo presidente e collaboratori. 
E intanto il libro fa ancora bella figura sul sito dell'ISTAT, oltre che a essere venduto online e in libreria. Brutta storia, direi\ldots\ E le autrici?
\item[G.] Le ho contattate, ovviamente --  è il modo con cui ci si comporta in ambiente scientifico  --  segnalando in modo puntuale quello che secondo me non andava nella fiaba delle streghe, a cui ero saltato rapidamente dopo aver letto le prime due, perché è il tema che mi stava più a cuore. Poi ho smesso perché mi si son drizzati i capelli  --  si fa per dire --  e mi son limitato a segnalare la fiaba a colleghi e amici esperti di probabilità e inferenza bayesiana. I quali
  sono rimasti inorriditi. Ma si son fatti anche grande risate,
  essendo il tipico caso in cui non sai se ridere o piangere:
  ti vien da ridere sui responsabili, ma da piangere sulle conseguenze. 
\item[N.] E loro, le autrici, ti hanno risposto?
\item[G.] Sì, almeno loro, sì. Una, per la precisazione,
  ma anche a nome delle altre, suppongo, visto che erano in copia.
\item[N.] E che dicevano, o che diceva?
\item[G.] Essenzialmente ``si arrampicavano sugli specchi per difendere l'indifendibile''. Questo è stato il giudizio unanime dei miei amici
esperti di queste cose, ai quali avevo inoltrato la risposta. E intanto il libro sta ancora sul sito dell'ISTAT, è raccomandato in siti vari dove le notizie vengono ripostate senza nemmeno che uno si prenda la briga di valutare o far valutare ciò che pubblicizza, e, soprattutto, credo circoli nelle scuole.
\item[N.] Sì, avevamo visto, cercando su internet. Addirittura abbiamo trovato un dirigente scolastico che diligentemente lo raccomandava a alunni e genitori con tanto di circolare protocollata. È per quello che ci aveva interessato.
\item[A.] Poi, fra quello che non capivo io e i dubbi di Noemi, abbiamo deciso di contattarti, visto che ogni tanto capiti da queste parti. 
Ti ringraziamo, Giulio, mi sembra che ci siamo chiarite. Soprattutto spero bene che Noemi si ricordi tutto, perché io non sarei in grado di ripetere molto di quanto ci siamo detti, anche se mi sembra di aver afferrato la sostanza. Poi vedremo cosa fare a scuola.
\item[N.] Scusate, ma ci siamo dimenticati una cosa, il famoso problema predittivo, inteso come valutazione della probabilità del colore del cappello che vedremo una certa mattina all'ingresso della caverna, a partire dalla sequenza dei colori osservati nei giorni precedenti.
\item[G.] Non me ne ero dimenticato, ma siccome l'abbiamo fatta un po' lunga e abbiamo riempito tutti questi fogli, non volevo esagerare. È sicuramente un problema interessante. Soprattutto per il legame con il problema inferenziale. Per quello parlavo di inferenziale-predittivo. Ma se volete\ldots\  non ci vuole molto, niente di complicato.
\item[N.] A me andrebbe, se non vi dispiace, ora che abbiamo fatto trenta\ldots
\item[A.] Mamma mia, qui mi ci vuole altra tisana\ldots\  credo per tutti.
  E porto pure degli aperitivi, vista l'ora.
\end{description}
%
\begin{center}
{\Large  $[\,$ Pausa $\,]$ }
\end{center}  
\mbox{}\\

\begin{description}
\setlength\itemsep{0.5mm}
\item[G.] Affrontiamo allora il problema predittivo, che nella nostra storia si riferisce al cappello che le streghe esporranno domani. Concentriamoci sul colore viola, V, come abbiamo fatto prima. Cosa ci aspettiamo il primo giorno?
\item[N.] Stessa probabilità di Viola e di Nero, direi.
\item[G.] E se nei primi 10 giorni abbiamo visto cinque volte V e cinque volte N, quanto varrà la probabilità di V l'undicesimo giorno?
\item[A.] 50\%, o c'è il trucco? Con te non si sa mai\ldots\  Simmetria iniziale, simmetria nelle osservazioni: rimaniamo in condizioni di simmetria, mi sembra.
\item[G.] Non lo potevi dire meglio, a parte precisare doverosamente 
che abbiamo un fattore di aggiornamento che dipende solo dalla differenza e quindi se la differenza è nulla il fattore di aggiornamento è unitario.
E se avessimo osservato 50 V e 50 N in cento giorni?
\item[A.] Stessa cosa, immagino. Dove vuoi andare a finire?
\item[G.] Da nessuna parte, almeno finché vediamo un ugual numero di cappelli dei due colori. Ma ora immagina di avere visto per 10 giorni di seguito il cappello viola, cosa ti aspetti l'undicesimo giorno?
\item[A.] Viola al 100\%? È qui che sta il trucco?
\item[G.] Fuoco\ldots\  Non perché sia il 100\%, ma perché hai sospettato il trucco, come dici tu. Ma si tratta di un trucco a fin di bene.
\item[N.] 100\% sicuramente non può essere. Un minimo di incertezza ci deve pur essere\ldots 
\item[G.] Va bene, escludiamo 100\%, per ovvi motivi. Ma se proponessi 99\% ti andrebbe bene? 
\item[N.] Uhm, non so, anche se mi sembra alta. 
\item[G.] Va bene, diciamo che 99\% è alta. 98\%? 95\%? 90\%, 85\%, 80\%\ldots
\item[A.] Se continui al ribasso, quando arrivi a qualche percento dirò che sicuramente è sbagliata. Non so perché, ma me lo sento, ah ah\ldots
\item[G.] Scusate, detti così sembrano numeri senza senso e non è chiaro qual è quello giusto, senza ragionarci un po' su. Ma ricominciamo da capo. Quali erano le composizioni possibili, ovvero il numero di streghe con il cappello viola, essendo l'altro numero il complementare per raggiungere 21?
\item[A.] Avevamo detto che potevano essere 7 o 14 viola, su un totale di 21, e avevamo chiamato le due possibilità $V_7$ e $V_{14}$.
\item[G.] E quanto vale la probabilità di cappello viola nei due casi?
\item[A.] L'avevamo calcolato prima. Anzi le aveva calcolato Noemi. Ecco, su quel foglio: 1/3 e 2/3. 
\item[G.] E quindi?
\item[N.] Capito!! Se fossimo sicuri di $V_7$ avremmo 1/3, ossia circa il 33\%. Se invece fossimo sicuri di $V_{14}$ avremmo 2/3, ossia circa il 67\%. Inizialmente eravamo ugualmente incerti sulle due possibilità, e avevamo quindi 1/2, o 50\%, a metà strada fra i due valori, anzi, più precisamente il valore medio.
\item[G.] FUOCO!!
\item[N.] Quindi è chiaro che la probabilità di viola non potrà mai essere superiore a 2/3, e quindi nessuno dei valori che ci avevi proposto prima era accettabile. Che scema! E per lo stesso motivo non potrà mai essere inferiore a 1/3, sempre assumendo le ipotesi del problema.
\item[G.] Inevitabilmente. E quindi? Mi sembra che ci sei arrivata. Basta solo che connetti questi ragionamenti a qualche formula che devi aver studiato. Anzi, forse se ci pensi un attimo non hai nemmeno bisogno di ricordare la formula. Prova.
\item[N.] Guarda, direi che più crediamo a $V_{14}$ e più la probabilità di osservare un cappello viola si avvicina a 2/3. Dal fattore di Bayes-Turing del problema abbiamo che a ogni osservazione di viola la probabilità di $V_{14}$ raddoppia rispetto a quella di $V_7$. Dopo dieci volte abbiamo $2^{10}$, che fa 1024, ovvero $V_{14}$ diventa circa mille volte più credibile di $V_7$. Quindi la probabilità di viola il giorno dopo diventa praticamente pari a $2/3$, appena un pochino meno.
\item[G.] E se volessi arrivare al numero preciso?
\item[N.] A intuito direi che farei una media fra $2/3$ e $1/3$, \ldots ma non una semplice media, bensì una media ponderata\ldots
\item[G.] Che brutto termine, usato sì nei libri di testo, ma sembra con l'intento di spaventare i ragazzi! Tu, quando vai al mercato, la frutta te la fai ponderare?
\item[N.] Ok, `media pesata', nella quale i pesi sono proporzionali a quanto credo a ciascuna delle due ipotesi.
\item[G.] Sì, ottimo! Hai espresso a parole l'essenza della formula della teoria della proba\-bilità alla quale mi riferivo, anche se in genere nei testi non viene presentata come media pesata, o ponderata che sia, peccato\ldots \\
Allora il numero `esatto'? Insomma scritto con un certo numero di cifre da poterne apprezzare la differenza rispetto a 0.66666666\ldots \\
Puoi usare una calcolatrice  --  sugli smart non mancano.
\item[N.] Preferisco i conti esatti, finché è possibile, e scrivo qui i passaggi.\\
Siccome
          $$\frac{P(V_{14}\,|\,10V) }{ P(V_7\,|\,10V)} = \frac{1024}{1}$$
abbiamo
          $$P(V_{14}\,|\,10V) = \frac{1024}{1025}$$
e
          $$P(V_7\,|\,10V) = \frac{1}{1025}$$
e, poiché la loro somma vale chiaramente uno, nella formula della media pesata nemmeno staremo a dividere per la somma dei pesi, in quanto unitaria. Abbiamo quindi
\begin{eqnarray*}
P(V\,|\,10V) &=& P(V\,|\,V_{14}) \times P(V_{14}\,|\,10V) + P(V\,|\,V_7) \times P(V_7\,|\,10V) \\
             &=& \frac{2}{3} \times \frac{1024}{1025} + \frac{1}{3} \times \frac{1}{1025} \\
               &=& \frac{2048}{3075} + \frac{1}{3075} \\
               &=& \frac{2049}{3075}
\end{eqnarray*}
e a questo punto mi serve la calcolatrice. Un attimo che prendo il telefonino. Ecco, 0.6663415\ldots, ma non esageriamo, mi fermo dove comincia a differire significativamente da 2/3:  0.6663. Quindi la probabilità di cappello
viola anche l'undicesimo giorno vale il 66.63\%, giusto?
\item[G.] L'hai calcolata tu e, come abbiamo detto, ci aspettavamo dovesse venire un 'nticchia sotto 2/3, come è risultato dai conti. E per concludere, ecco la
  rete bayesiana delle questioni inferenziale-predittive legate alla fiaba delle
  streghe.
\begin{center}
\epsfig{file=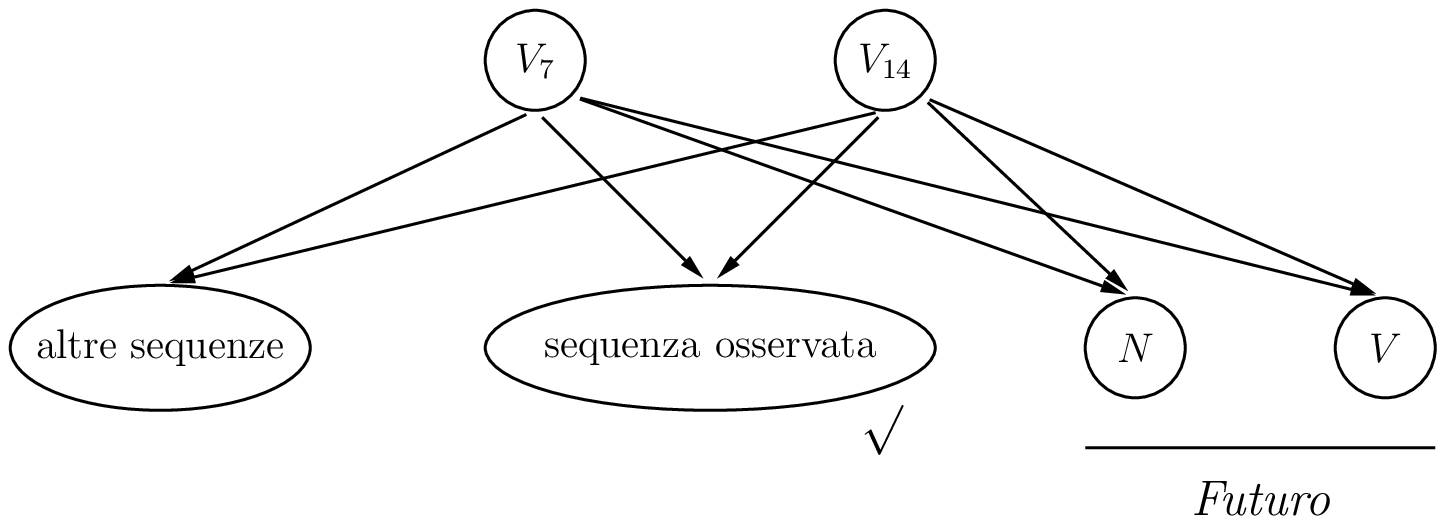,clip=,width=\linewidth}
\end{center}
A partire da una sequenza osservata, con simbolo di spunta $\surd$ a ricordare
che è `accertata', inferiamo i parametri liberi del modello,
che sono appunto le probabilità delle ipotesi $V_7$ e $V_{14}$,
e da queste effettuiamo previsioni probabilistiche su eventi futuri.
\item[A.] Sempre più intriganti queste reti bayesiane.  E quindi
  immagino che possiamo pure effettuare una previsione probabilistica
  sul cibo che si aspetta la strega. Basta aggiungere altri
  due tondi con le rispettive frecce. Posso? Ecco, avendo cancellato
  la tua scritta sotto, che non ci serve più.
\begin{center}
\epsfig{file=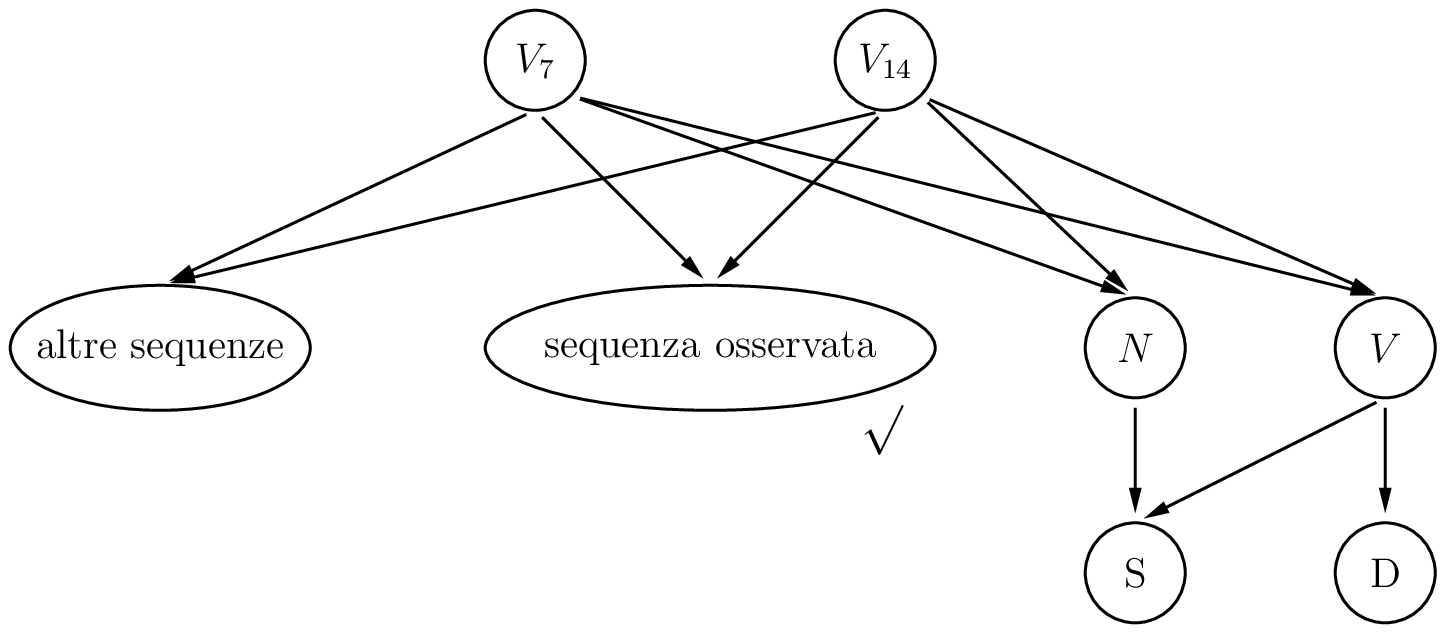,clip=,width=\linewidth}
\end{center}
E da questa uno dovrebbe essere in grado di calcolare
la probabilità che la strega estratta a sorte preferisca
Dolce o Salato, giusto.
\item[G.] E nel caso di una sequenza di dieci viola non è nemmeno
  troppo difficile.
\item[N.] Sì, perché siamo praticamente certi che si tratti di $V_{14}$,
  quindi le probabilità di $N$ e di $V$ sono sostanzialmente
  1/3 e 2/3. Ne segue che la probabilità che la strega estratta preferisca
  Dolce à data da
  $$\frac{6}{7} \times \frac{2}{3} = \frac{12}{21} = \frac{4}{7}\,,$$
  mentre la probabilità che preferisca Salato vale
  $$1 \times \frac{1}{3} + \frac{1}{7}\times \frac{2}{3}$$
  ovvero 
  $$ \frac{1}{3} +  \frac{2}{21} =  \frac{7+2}{21} =  \frac{9}{21} = \frac{3}{7}\,.$$
  E la somma dà giustamente uno perché o preferisce Dolce o Salato.
\item[G.] Ottimo!
\item[A.] Scusate, ma se leggo bene quelle frazioni,
  stiamo dicendo che delle streghe nella grotta a dodici su ventuno piace dolce
  e a nove su ventuno piace salato.
  Era quello che avevo detto io prima senza tutti quei conti!
\item[G.] E la cosa ci conforta, in quanto il tuo ragionamento
  era valido se assumevi quattordici streghe dal cappello Viola
  e sette dal cappello Nero. E adesso, come abbiamo detto, siamo
  praticamente certi, e in effetti abbiamo usato tale approssimazione. 
  I conti che sembrano complicati
  servono quando siamo in stato
  di incertezza.
  Prendiamo il caso della sequenza di quattro neri consecutivi che
  abbiamo visto prima, che ovviamente ci faceva protendere a favorire
  $V_7$, ma con un margine di incertezza più elevato.
\item[N.] Certo le probabilità di $V_7$ e di $V_{14}$ erano
  rispettivamente 16/17 e 1/17. E se prendiamo il caso
  di Salato, dove ci sono tutti i termini, abbiamo  
  \begin{eqnarray*}
    P(S\,|\,4N) &=& P(S\,|\,N)\cdot P(N\,|\,4N) +
                    P(S\,|\,V)\cdot P(V\,|\,4N) \\
                    &=& P(S\,|\,N)\cdot
                    \left[ P(N\,|\,V_7)\cdot P(V_7\,|\,4N) +
                      P(N\,|\,V_{14})\cdot P(V_{14}\,|\,4N)\right] + \\
                    &&  P(S\,|\,V)\cdot
                    \left[ P(V\,|\,V_7)\cdot P(V_7\,|\,4N) +
                      P(V\,|\,V_{14})\cdot P(V_{14}\,|\,4N)\right] \\
                    &=& 1\times  \left[ \frac{2}{3} \times \frac{16}{17}
                      + \frac{1}{3} \times \frac{1}{17}\right]
                    \ + \ \, \frac{6}{7}\times \left[ \frac{1}{3} \times \frac{1}{17}
                      + \frac{2}{3} \times \frac{16}{17}\right]\\
                    &=& \ldots 
  \end{eqnarray*}
  \mbox{}\\ \vspace{-1.4cm} \mbox{}
\item[A.] Va bene, basta così, mi arrendo!
\item[N.] Un altro problema interessante potrebbe essere
  quello di calcolare la probabilità del colore esposto,
  se veniamo a sapere soltanto che la strega ha gradito
  un particolare cibo.
\item[G.] Perché no? Come si dice,
  l'appetito vien mangiando, visto che si parla di cibo.
  Anzi, vi propongo di provare
  a risolverlo assumendo  che il cibo gradito sia stato
  Salato e sotto due ipotesi: la prima è quando
  le ipotesi $V_7$ e $V_{14}$ ci sembravano equiprobabili;
  la seconda è
  quando siamo certi si tratta di $V_{14}$. E mi raccomando di fare
  anche il diagramma.
\item[N.] Certo, ci proveremo.
\item[A.] Sì, va be', ci proverai tu\ldots\  io posso aiutare a 
  tracciare la rete bayesiana, visto che ci ho preso la mano. 
\item[N.] Vedrai che non sarà difficile.
  Ma, scusate, ora sono io ad avere ancora un dubbio.
  Riguarda la differenza fra la probabilità di Viola
  il giorno dopo fatto in questo modo e quello che avrei fatto a partire dalla
  frequenza relativa di tale colore nel passato,  
  che in questo caso mi sembra essere addirittura sbagliata.
  Ma forse è immagino perché il numero di osservazioni nel passato è piccolo, solo dieci.
\item[G.] Anche se fosse stato cento, capisci bene che non è improbabile avere frequenze di viola, intendo quelle osservate, sotto il 33\% o sopra il 67\%, dando luogo a valori non accettabili per la probabilità di cappello viola alla luce di quanto sappiamo.
\item[N.] Ma allora, la valutazione di probabilità basata sulla frequenza di successi verificatesi nel passato?
\item[G.] Dipende, ma certamente non in questo caso, stanti le nostre ipotesi. Non mi sembra ci siano dubbi. Mentre in altri casi può andare, e si basa proprio su un ragionamento bayesiano, ma purtroppo non possiamo scendere nei dettagli. Dico soltanto che andrebbe bene se avessimo un grandissimo numero di streghe, con la proporzione di quelle con il cappello viola che possa
  andare `praticamente con continuità' da zero al 100\%. E la formula che dà la probabilità del prossimo cappello viola, avendo osservato $x$ cappelli viola in $n$ giorni, vale $(x+1)/(n+2)$, approssimabile a $x/n$ in tanti casi pratici. Solo in questo caso, e con queste condizioni, la probabilità di viola il giorno successivo sarebbe pari alla frazione di viola osservata nel passato. E, per inciso, $(x+1)/(n+2)$ è nota come {\em formula di successione di Laplace}.
  -- Con questa nota storica, che fra l'altro ci giustifica l'uso delle frequenza per valutare la probabilità, direi che abbiamo concluso.
\item[N.] Ma una cosa ce la puoi ancora dire, forse: come si raccorda la formula di Bayes che si trova in giro con quella `tua', per così dire, di aggiornamento dei rapporti di probabilità.
\item[G.] Non c'è problema e ci vuole veramente poco. Ecco, cerchiamoci su internet una formula. Suggerirei di andare sul sito dell'ISBA, la società internazionale di analisi
  bayesiana\cite{ISBA} --  chi meglio di loro!
Eccola, al neon, 
\begin{center}
\epsfig{file=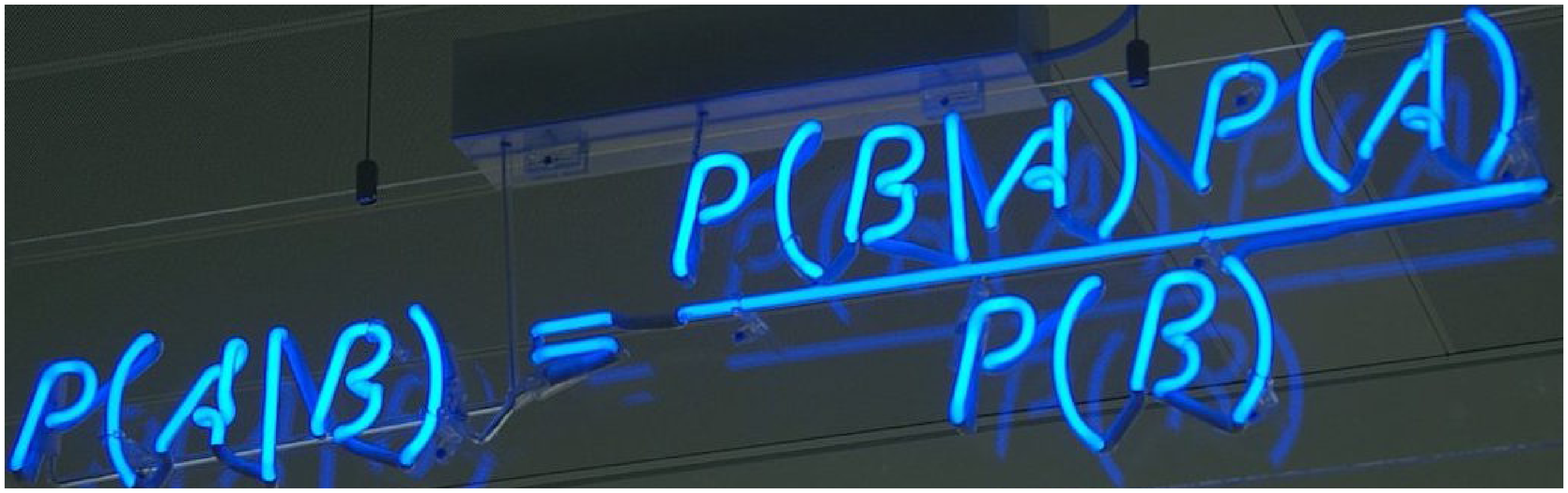,clip=,width=\linewidth}
\end{center}
\item[A.] Ma questa l'ho vista! Non stava su una
  maglietta che portavi qualche estate fa? 
\item[G.] E già, me l'avevano regalata degli amici che si occupano in modo professionale di inferenza e previsioni probabilistiche.
\item[A.] E non sta pure sulla copertina di un tuo libro? Quello del
  Mago\ldots\  ma non di Oz, \ldots\ era qualcosa del genere\ldots 
\item[G.] Sì, qualcosa del genere, di Odds, che suona quasi uguale\,\cite{MagoOdds}.
  Gli {\em odds} in inglese sono i rapporti di probabilità, quelli di cui abbiamo parlato -- per la precisione di una ipotesi contro l'ipotesi alternativa.
\item[A.] Ah, ecco tante cose che si incastrano insieme.
\item[N.] È la stessa che avevo trovato io su
  internet,\footnote{\,\url{https://learningtogo.info/2017/11/16/your-bayesian-brain-are-you-wired-for-statistics/}} con sullo sfondo un ritratto
  del reverendo Thomas Bayes e una citazione che dà una idea, seppur molto generale,
  del ragionamento che porta il suo nome. Eccola, ce l'ho sul telefonino, guardate: 
\begin{center}
\epsfig{file=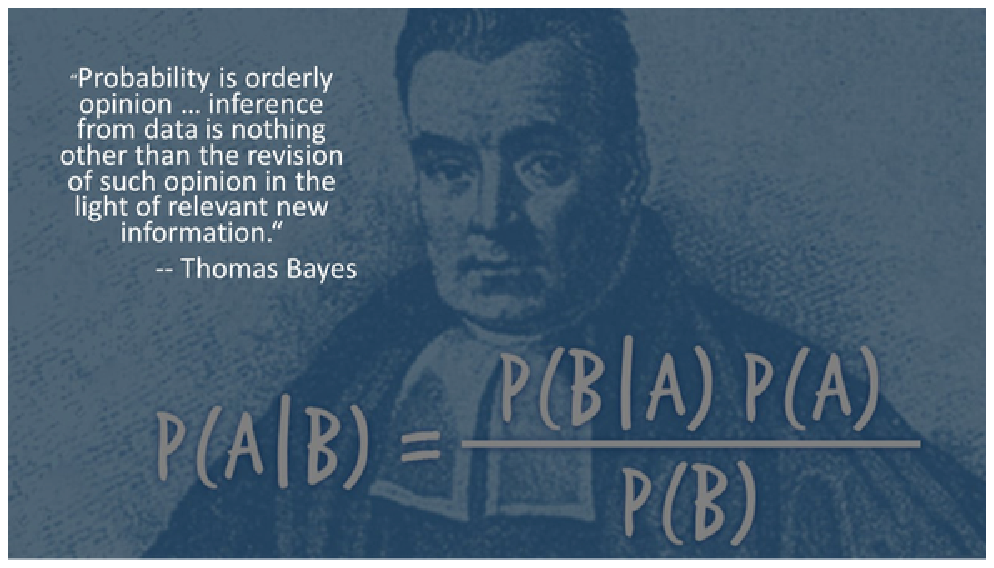,clip=,width=\linewidth}
\end{center}
\item[A.] Bello il concetto di ``opinione ordinata''!
  Non semplicemente vero o falso,
  ma con delle gradazioni intermedie.
  Che è quello che intendiamo sempre quando diciamo che un'ipotesi
  è più o meno probabile. 
  E come si passa dalla formula `al neon' a quella `del Mago', in cui entrano {\em odds} e fattore di Bayes-Turing?
\item[G.] Facile. Innanzitutto, esprimendoci in termini di cause e effetti, $B$ sta per l'effetto osservato, ``{\em the relevant new
  information}'' della citazione. 
  $A$ sta invece per una causa di quell'effetto.
  $P(A)$ a destra è la probabilità iniziale di $A$, e quindi, con tutta la stima e il rispetto per gli amici dell'ISBA, la scriverei $P_0(A)$.
  Invece $P(A\,|\,B)$ a sinistra ci dice come tale probabilità è aggiornata da $B$.
\item[A.] Ma non vedo altre cause in quella formula!
\item[G.] Ci sono, ci sono, ma sono nascoste nella $P(B)$ al denominatore, la quale va calcolata tenendo conto delle altre possibili cause. Se, ad esempio, oltre a $A$ abbiamo anche la causa alternativa $\overline A$, usando quella formula della media pesata che
  Noemi si era ricordata prima\ldots
\item[N.] La formula della probabilità totale\ldots 
\item[G.] Sì, un nome del genere, i nomi non sono il mio forte\ldots\  Dicevo se abbiamo due cause possibili, la probabilità dell'effetto $B$ la possiamo scrivere come
  $$ P(B) = P(B\,|\,A) \cdot P(A) + P(B\,|\,\overline A) \cdot P(\overline A)\,. $$
\item[A.] Ah, ho capito, con $\overline A$ che sta per non-$A$.
          {\em Tertium non datur}, come si diceva.   
\item[N.] Comunque, qualsiasi cosa sia, è chiaro che $P(B)$ al denominatore
  è lo stesso sia nel caso di $P(A\,|\,B)$ che in quello di
  $P(\overline A\,|\,B)$.
\item[G.] E quindi possiamo vedere $1/P(B)$ nella formula `al neon' come una costante moltiplicativa. Possiamo allora leggere la formula nel seguente modo:
  {\em la probabilità di $A$ dato $B$ è proporzionale alla probabilità iniziale di $A$ e alla probabilità di $B$ dato $A$}, anche se tale lettura è ancora di tipo matematico. Un modo migliore per capire l'essenza del ragionamento bayesiano, quando lo schema cause-effetti è applicabile, è che {\em la probabilità della causa $A$, alla luce dell'effetto $B$, è proporzionale non solo alla probabilità che tale causa possa dare un effetto del genere ma anche a quanto ritenevi probabile tale causa `a priori', ovvero indipendentemente dall'osservazione di tale effetto.}\footnote{\,Questa frase esprime in sostanza quello che
    Laplace chiama 
    nel capitolo 
    {\em Prin\-cipes généraux du Calcul des Probabilités}
    di Ref.~\cite{Laplace} 
    {``\em VI}\,°\,{\em Principio''}, che eleva a
    {\em ``principe fondamental de cette branche de l’analyse des hasards, qui consiste à remonter des évènemens aux causes''}\,({\bf !}). In
    formule 
    $$P(C_i\,|\,E) = \frac{P(E\,|\,C_i)\cdot  P_0(C_i) }{\sum_{i=1}^nP(E\,|\,C_i)\cdot P_0(C_i)}
    = \frac{P(E\,|\,C_i)\cdot  P_0(C_i)}{P(E)}\,$$
    ove $P(E) = \sum_{i=1}^nP(E\,|\,C_i)\cdot P_0(C_i)$ è una estensione
    di
    $ P(B) = P(B\,|\,A) \cdot P(A) + P(B\,|\,\overline A) \cdot P(\overline A)$
    di p.\,43 
    a $n$ cause $C_i$ mutuamente esclusive. {\bf Da cui segue},
    essendo $P(E)$ una costante numerica, dati il modello, espresso
    nelle varie $P(E\,|\,C_i)$, e le prior $P_0(C_i)$,
    $\mathbf{ P(C_i\,|\,E) \propto P(E\,|\,C_i)\cdot  P_0(C_i)}$.\\
    (Vedi anche \url{https://www.youtube.com/watch?v=YrsP-h2uVU4&t=146s} da 1:15:10.)
%
  }
\item[A.] Interessante questa lettura. Un esempio?
\item[G.] Pensate a sintomi di possibili malattie: una molto comune, che potrebbe dare quei sintomi con un una certa probabilità; l'altra molto rara, ma che provoca sempre quei sintomi. Qualsiasi medico minimamente esperto sa che deve pensare subito alla malattia più comune. Insomma, come diceva non ricordo chi in uno dei primi episodi  di {\em Dr. House}, ``Quando senti rumore di zoccoli, pensa al cavallo, non alla zebra''.\footnote{\,Per l'attribuzione del detto
  vedi \url{https://it.wikiquote.org/wiki/Theodore_Woodward}.} 
\item[A.] Giusto, se stai in America o dalle nostre parti! Anzi, dalle nostre parti, una cinquantina di anni fa si sarebbe detto ``pensa all'asino''.
\item[G.] Ah ah, giusto! Ma, comunque, a pensarci bene, 
che significa rumore di zoccoli? Dipende dalla stazza dell'animale
e dal suo incedere, e se è ferrato o no. Io sulle nostre strade 
di campagna riconoscevo facilmente se si trattava di asino o di cavallo,
al più potevo essere incerto fra cavallo e mulo.
Quindi il detto americano lo vedo valido solo in qualche prateria, 
o comunque su sterrato.\\
Tornando all'inferenza delle malattie dai sintomi, il fatto che la malattia rara dia con certezza sintomi del genere non implica che tali sintomi, osservati su un paziente, derivino con certezza da essa.
\item[N.] Fortunatamente, direi! Anche perché la malattia rara è in genere associata a malattia grave. E tornando alla parte formale, mi sembra di aver capito come si passa della formula `al neon' 
  a quella con il fattore di Bayes-Turing.
\item[A.] Beata te. 
\item[N.] È facile. Immagina una formula analoga, nella quale,
  invece di comparire $A$, compaia l'ipotesi alternativa $\overline A$:
  $$P(\overline A\,|\,B) = \frac{P(B\,|\,\overline A) \cdot
    P(\overline A)}{ P(B)}\,.$$
  Basta dividere le due equazioni membro a membro e otteniamo
  $$ \frac{P(A\,|\,B)}{P(\overline A\,|\,B)} =
  \frac{\frac{P(B\,|\,A) \cdot  P(A)}{ P(B)}
       }
       {\frac{P(B\,|\,\overline A) \cdot  P(\overline A)}{ P(B)}
       }
       = \frac{ P(B\,|\,A) \cdot  P(A)
       }
       {P(B\,|\,\overline A) \cdot  P(\overline A)
       }\,,
       $$
       che possiamo riscrivere come
       $$
       \frac{P(A\,|\,B)}{P(\overline A\,|\,B)} =
       \frac{P(B\,|\,A)}{P(B\,|\,\overline A)} \times
       \frac{P(A)}{P(\overline A)}  \,.
       $$
%
In modo compatto,  e usando l'espressione
Odds per i rapporti di probabilità e l'acronimo BTF per il fattore di Bayes-Turing,
otteniamo quindi la regola di aggiornamento che ci hai fatto vedere prima:
      $$   {\large \mbox{Odds}_{finali}  =  \mbox{BTF} \times \mbox{Odds}_{iniziali} } $$
\item[G.] Ben detto, quindi a questo punto veramente abbiamo ben poco da aggiungere sul ragionamento bayesiano, almeno per quanto riguarda la sostanza.
\item[A.] Grazie, Giulio. Ma chissà quanto ci ricorderemo di tutto questo. Io certo non molto, ma se mi aiuta Noemi a riordinare le idee, prometto che mi ci metto d'impegno.
\item[N.] Ma certo, e sarebbe simpatico se potessimo avere qualche esercizio da svolgere. In fondo siamo tutti insegnanti e sappiamo bene che solo quando uno prova a risolvere un problema si rende conto se veramente ha capito.
\item[G.] Sì, hai ragione e 
  uno ve lo avevo già proposto prima, in realtà su tua idea: quello
  di inferire il colore del cappello esposto sapendo che la strega ha gradito
  il cibo.  
  Ne posso suggerire uno al volo, di carattere medico, ma questa 
volta piacevole: niente malattie, bensì lieti eventi.
\item[A.] Ottima idea! E magari ci provo pure io perché effettivamente la matematica sembra semplice, a sapere cosa calcolare\ldots 
\item[G.] Eccolo, e scrivetevi pure i dati numerici, così è più facile controllare se lo avete fatto bene.
  Immaginate un test per conoscere, seppur con incertezza, il sesso del nascituro. Ovviamente non si tratta di ecografia, ma di qualche analisi clinica che si possa effettuare anche quando con l'ecografia non si capirebbe niente. Ecco i dati
  \begin{itemize}
    \item   se il nascituro è {\bf M}aschio il risultato del test dà {\bf m}aschio al 95\% 
       (e al 5\% {\bf f}emmina); 
    \item   se il nascituro è {\bf F}emmina il risultato del test dà {\bf f}emmina all'80\% 
       (e al 20\% {\bf m}aschio); 
    \item   per semplicità assumiamo che il feto abbia la stessa probabilità
    di essere {\bf m}aschio o {\bf f}emmina.
    \end{itemize}
Suggerirei quindi di usare come variabili del problema $M$, $m$, $F$ e $f$, con ovvio significato.
Immaginate ora che due signore, Anna e Barbara, eseguano tale test. Il risultato di Anna è `maschio' e quello di Barbara è `femmina'. Quale delle due signore dovrà essere più confidente sul sesso che le è stato anticipato?
\item[A.] Direi quella a cui è stato detto maschio, {\em ovviamente}, perché --  ecco, fatemelo scrivere in formule,  $P(m\,|\,M) = 95\%$, mentre $P(f\,|\,F) = 80\%$. 
\item[G.] Né ovvio né tanto meno corretto. Come vedete, su problemi di questo è facile sbagliarsi. E alla grande!
\item[A.] Non è così? Me lo sentivo che con te non c'è da fidarsi e che avrei dovuto dire
  il contrario di quello che pensavo\ldots\ 
   Ma allora non ho capito niente!
\item[G.] Secondo me il ragionamento bayesiano lo hai capito. Quello che non hai ancora capito è che devi usarlo, con calma.  
  Provate a far i conti e poi ne riparliamo. Vedrete che sono facili. E magari fateli in entrambi i modi, sia usando i rapporti di probabilità che usando la `formula al neon'  --  chi è che l'aveva chiamata così? A saperlo venivo con la maglietta bayesiana, ma sotto camicia e pullover, ovviamente, data la stagione.
\item[N.] Sì, dovremmo aver capito, speriamo, e ci possiamo mettere alla prova facendo l'esercizio che ci hai lasciato. E magari ci puoi suggerire qualcosa da leggere. Il problema è cosa fare a scuola. Veramente in questo momento mi sento confusa, più su questo che sul ragionamento bayesiano in sé.
\item[G.] Sì, lo capisco e non ho grandi suggerimenti da dare, o, almeno, consigli operativi su come procedere. Dovete provare voi, per capire quali sono le difficoltà che incontrerebbero i ragazzi. Poi ne riparliamo. Capisco bene che tante cose che a me sembrano scontate in realtà non lo sono affatto.
Già scrivere formule con dei simboli può essere per tanti un ostacolo insormontabile. Ma conoscendo il {\em target}, come si dice oggi, con cui avete a che fare, forse vi viene qualche idea.
\item[A.] Oppure potremmo, tanto per cominciare, fare la sperimentazione didattica a cui avevamo pensato prima di parlare con te. Pensavamo infatti cominciare
  proponendo la lettura delle fiabe a gruppetti
  e poi valutare cosa hanno capito.
\item[N.] E per allargare il campione potremmo  anche coinvolgere altri
  colleghi. Potrebbe essere interessante.
\item[A.] Comunque ci hai convinte che nella fiaba delle streghe
  il ragionamento bayesiano non c'è. A questo si aggiunge il fatto che la strategia decisionale è sbagliata, e non so che altro che non andava. E anche dal punto narrativo la storia zoppica un po'. Quindi questa fiaba nemmeno la assegnerei nei test di comprensione.
\item[G.] Ma forse potreste cambiarla, o inventarne un'altra dalla quale si possa evincere la sostanza dell'aggiornamento delle probabilità.
\item[N.] Ci possiamo provare, ma non mi sembra facile, almeno in questo momento. Chissà se parlandone e interagendo con ragazzi e colleghi ci verranno delle idee.
\item[G.] Sì, io pure in questo momento sono un po' scettico, e comunque sono dell'avviso di non forzare la mano. Meglio non insegnare che insegnare cose sbagliate.
\item[A.] Sì, lo so che sei contrario alla divulgazione, ne abbiamo parlato.
\item[N.] Sul serio? Questa non la sapevo. Ma ora, parlando con noi, non stai facendo divulgazione, in un certo senso?
\item[G.] Sì e no, dipende da cosa si intende. Diciamo che sono per l'insegnamento e non per la divulgazione {\em cheap}, come direbbero gli inglesi, che poi è la quasi totalità di quella che circola, e di cui non vedo altro beneficio che far fare soldi a chi la vende. A questo proposito vorrei leggervi un paio di citazioni che esprimono bene il mio punto di vista.
  Le potete trovare in fondo alla mia pagina web dedicata
  all'insegnamento.\footnote{\,\url{http://www.roma1.infn.it/~dagos/teaching.html\#DEA}}
  Un momento che mi collego. Eccole. La prima è tratta da un pezzo di Giorgio Gaber di una ventina di anni fa, intitolato l'{\em Ingenuo}\,\cite{Gaber}:
  \begin{quote}
    {\sl 
     Tu credi che se un uomo ha un'idea nuova, geniale, \\
     abbia anche il dovere di divulgarla. Tu sei un ingenuo.\\
     Prima di tutto perché credi ancora alle idee geniali.\\
     Ma quello che è più grave è che tu credi all'effetto\\
     benefico dell'espansione della cultura.\\
     \ldots \   \ldots \\
     Tu mi dirai che la divulgazione è un dovere civile e che\\
     evolve il livello culturale della gente. Non riesci proprio\\
     a staccarti da un residuo populista e un po' patetico.\\
     Purtroppo, oggi, appena un'idea esce da una stanza è subito\\
     merce, merce di scambio, roba da supermercato. La gente se\\
     la trova lì, senza fatica, e se la spalma sul pane, come la Nutella.
     }
     \end{quote}
  Per la seconda devo ringraziare Alba, che mi ha fatto scoprire non molto tempo fa le {\em Lettere di Berlicche} di Clive Lewis\,\cite{Berlicche}, quello delle più famose {\em Cronache di Narnia}. Si tratta di lettere che un vecchio e esperto diavolo scrive al diavoletto suo nipote, Berlicche appunto, per
  dargli consigli su come portare dalla sua parte la persona che gli è stata affidata. Ecco il finale della prima lettera, dalla quale mi sembra di evincere che la divulgazione sia opera di Satana.
 \begin{quote}
    {\sl 
      L'ideale è, naturalmente, di non fargli leggere neppure una riga di veramente scientifico, ma di infondergli l'idea generale grandiosa che egli conosce tutta la scienza, e che ogni cosa gli avvenga di raccogliere in conversazioni casuali e nelle lettura è ``i risultati della moderna investigazione'', Ricordati bene che il tuo dovere è di ubriacarlo. Dal modo con il quale alcuni di voi giovani demoni parlate si potrebbe pensare che la vostra occupazione sia quella di insegnare. {\footnotesize \\
      \mbox{} \\
      \mbox{} }
      \hspace{0.5cm} Tuo affezionatissimo zio.
       {\footnotesize \\
      \mbox{} }
    }
    \end{quote}
  
\end{description}


\end{document}